\documentclass{article}

\usepackage{natbib} 

\usepackage[english]{babel}
\usepackage[utf8]{inputenc} 
\usepackage[T1]{fontenc} 

\usepackage{amssymb} 
\usepackage{amsthm} 
\usepackage{amsmath} 
\usepackage{mathtools} 
\usepackage{dsfont} 
\usepackage{mathrsfs} 
\usepackage{stmaryrd} 
\numberwithin{equation}{section} 

\usepackage[nottoc]{tocbibind} 
\usepackage{hyperref} 

\usepackage{enumitem} 
\usepackage{array} 

\usepackage{geometry} 

\allowdisplaybreaks 

\newcommand{\noi}{\noindent} 

\newcommand{\prth}[1]{\left(#1 \right)} 
\newcommand{\crch}[1]{\left[#1 \right]} 
\newcommand{\ens}[1]{\left\{ #1 \right\}} 

\newcommand{\gO}[2]{\mathcal{O}_{#2}\left(#1 \right)} 

\newcommand{\scal}[2]{\left\langle  #1 , #2 \right\rangle} 

\newcommand{\tr}{\mathrm{Tr}} 
\newcommand{\supp}{\mathrm{supp}} 

\newcommand{\diff}{\mathrm{d}} 

\newcommand{\matp}[1]{\begin{pmatrix}#1\end{pmatrix}} 

\newcommand{\norme}[1]{\left\lVert #1\right\rVert} 
\newcommand{\normi}[1]{\left| \left| #1 \right| \right|_\infty} 

\newcommand{\abs}[1]{\left|#1\right|} 

\newcommand{\Sp}{\mathrm{Sp}} 
\newcommand{\sgn}{\mathrm{sgn}} 
\renewcommand{\L}{\mathrm{L}^2} 
\renewcommand{\H}{\mathrm{H}_\ep} 
\newcommand{\Op}{\mathrm{Op}_\ep} 
\newcommand{\Opp}{\mathbf{op}^\mathcal{C}_\ep} 
\newcommand{\Opd}{\mathrm{Op}_{\sqrt\ep,1}} 
\newcommand{\ep}{\varepsilon} 
\newcommand{\E}{\mathrm E} 
\newcommand{\D}{\mathrm{D}} 
\newcommand{\m}{\mathrm{m}} 
\renewcommand{\i}{\mathrm{i}} 

\newcommand{\N}{\mathbb{N}} 

\newcommand{\Z}{\mathbb{Z}} 

\newcommand{\R}{\mathbb{R}} 
\newcommand{\C}{\mathbb{C}} 
\renewcommand{\S}{\mathbb{S}} 

\newtheorem{theorem}{Theorem}[section]
\newtheorem{lemma}[theorem]{Lemma}
\newtheorem{corollary}[theorem]{Corollary}
\newtheorem{proposition}[theorem]{Proposition}

\theoremstyle{definition}
\newtheorem{remark}[theorem]{Remark}
\newtheorem*{rem*}{Remark}

\newtheorem{example}[theorem]{Example}
\newtheorem*{exe*}{Example}

\newtheorem*{def*}{Definition}
\newtheorem{assumption}[theorem]{Assumption}

\usepackage{sectsty}
\sectionfont{\centering}

\usepackage[rightcaption]{sidecap}
\sidecaptionvpos{figure}{c}
\usepackage{graphicx} 

\title{Propagation of Semiclassical Measures \\ Between Two Topological Insulators}
\author{\'ERIC VACELET}

\begin{document}

\maketitle

\begin{abstract}
\noi We study propagation in a system consisting of two topological insulators without a magnetic field, whose interface is a non-compact, smooth, and connected curve without boundary.
The dynamics are governed by an adiabatic modulation of a Dirac operator with a smooth, effective variable mass.
We determine the evolution of the semiclassical measure of the solution using a two-scale Wigner measure method, after reducing the Hamiltonian to a normal form.
\end{abstract}

\noi \textit{Keywords}: propagation, Dirac operator, semiclassical analysis, Wigner measures, semiclassical measures valued in a separable Hilbert space, time-dependent two-scale semiclassical measures.

\noi $2024$ \textit{Mathematics Subject Classification} : 35B40, 35F05, 35Q40.

\tableofcontents

\section{Introduction}

\paragraph{Topological insulators and Dirac-type equations.}
Topological insulators are electronic materials that possess a bulk band gap, like ordinary insulators, but exhibit protected conducting states at their edges.
The Dirac equation arises naturally in the context of topological insulators for honeycomb structures~\cite{FLW16, Fefferman_1, Drouot_topoinsu}, such as graphene~\cite{RH08, LWZ19, DW20}.

\medskip

In this paper, we study the following system of Dirac evolution equation on $\L\prth{\R^2,\C^2}$,
\begin{equation} \label{eq:Dirac}
\left\{
\begin{array}{ll}
	\left( \ep \D_t + \begin{pmatrix}  \m(x) & \ep \D_1 - \i\ep \D_2 \\ \ep \D_1 + \i\ep \D_2 & - \m(x) \end{pmatrix} \right) \psi^\ep  (t,x) = 0, & \forall (t,x) \in \R \times \R^2, \\
	\psi^\ep(0, x) = \psi^\ep_0(x), & \forall x \in \R^2,
\end{array}
\right.
\end{equation}
where $\displaystyle{\D_\# \coloneqq -\i \partial_\#}$ for $\# \in \{t,1,2\}$ with $\partial_j \coloneqq \partial_{x_j}$ for $j \in \{1,2\}$, and where $\displaystyle{\ep}$ is a small positive parameter.
The family $\displaystyle{\prth{\psi^\ep_0}_{\ep > 0} }$ is uniformly bounded (with respect to $\ep$) in $\displaystyle{ \L\prth{\R^2, \C^2} }$.
The function $x \mapsto \m(x)$ is a smooth real-valued function, that is, $\m \in \mathscr{C}^\infty\prth{\mathbb{R}^2, \mathbb{R}}$, whose derivatives are all bounded.
The interface between the two topological insulators is defined by
\[
    \E \coloneqq \left\{ x \in \mathbb{R}^2 \mid  \m(x) = 0\right\}.
\]
To ensure that $\m$ separates $\mathbb{R}^2$ into two distinct regions, we assume that $\E$ is a smooth, connected curve without boundary.
Under these conditions, there exists a unique family $\displaystyle{\prth{\psi^\ep_t}_{\ep > 0} }$ of solutions in $\displaystyle{ \mathscr{C} \left( \R,\L\prth{\R^2, \C^2} \right) }$ to~\eqref{eq:Dirac}.

\paragraph{Physical interpretation.}
The conducting states at the edge of topological insulators are called edge states.
Their propagation is one of the most important problems in the physical applications of topological insulators and superconductors~\cite{Ber13, Vol89}, photonics~\cite{LJS14, RH08, R13}, acoustics~\cite{PBSM15, MXC19}, and fluid mechanics~\cite{GJT21}.
These states are physically realized through the combined effects of spin–-orbit interactions and time-reversal symmetry.
More precisely, two-dimensional topological insulators are quantum spin Hall insulators, which are closely related to the integer quantum Hall state~\cite{HasanKane}.
In such materials, there exists an energy gap between the valence and conduction bands; however, unlike in a trivial insulator, these bands are ''twisted'' or ''crossing''~\cite{ZhuCheng}.
From a global point of view, Dirac-type equations often provide a simpler continuum macroscopic description of transport in a narrow energy band near a band crossing~\cite{Ber13, FC13, Vol89}.
Such models typically describe the evolution of the electronic wave function in topological insulators, such as graphene monolayers subjected to an external field or to non-homogeneous deformations.

\bigskip

Here, the topological insulators under consideration are characterized by the sign of $\m$: the region where $\m$ is positive corresponds to a topological insulator of index $1$, while the region where $\m$ is negative corresponds to a topological insulator of index $-1$.
The function $\psi^\ep(t,\cdot)$, also denoted by $\psi^\ep_t$, is interpreted as the wave function of an electron at time $t$ in a material consisting of two topological insulators.

\begin{SCfigure}[0.5][h]
\caption{Illustration of asymmetric transport for the mass $\m(x) = \arctan(x_1) - x_2$.}
\includegraphics[width=0.6\textwidth]{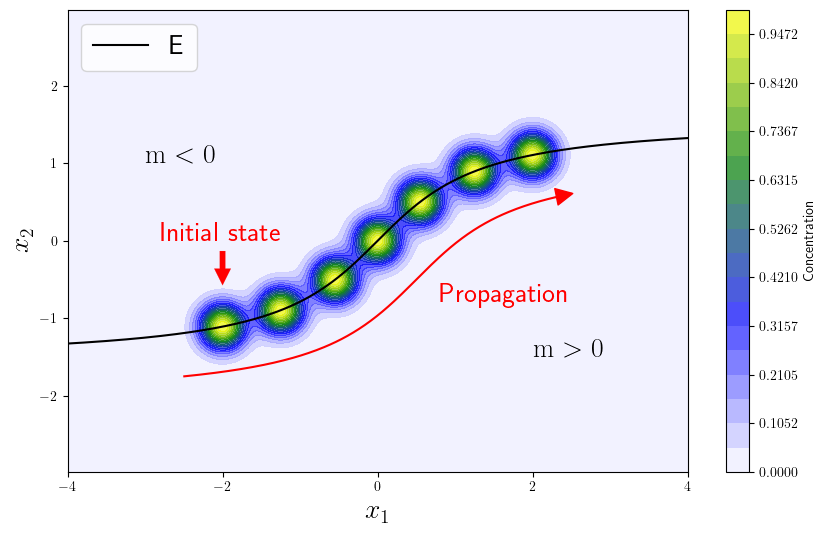}
\end{SCfigure}

Since equation~\eqref{eq:Dirac} originates from quantum mechanics through the quantum spin Hall effect, we are naturally led to study the problem in phase space $\mathrm{T}^*\R^2 = \R^2_x \times \R^2_\xi$.
In this setting, $\xi \in \R^2$ is interpreted as the momentum of the particle.

\subsection{Edge states generated by coherent states}

\paragraph{Coherent states.}
The system~\eqref{eq:Dirac} has already been studied in~\cite{Edge-States}, where the computed edge mode is generated by a very specific semiclassical wave packet.
Following~\cite{Drouot}, we define a wave packet concentrated at $(x_0,\xi_0) \in \R^4$, a point in phase space, as an $\ep$-dependent function satisfying, for all $\ep > 0$,
\[
    \mathrm{WP}^\ep_{x_0,\xi_0}[\vec{f}](x) \coloneqq \frac{e^{\frac{\i}{\ep} \xi_0\cdot (x-x_0)}}{\sqrt\ep} \vec f \prth{\frac{x - x_0}{\sqrt\ep}}, \quad x \in \R^2,
\]
with $\vec f \in \mathscr S \prth{\R^2,\C^2}$ independent of $\ep$.
We say that $\mathrm{WP}^\ep_{x_0,\xi_0}[\vec f]$ is oriented along $\vec V\in \C^2$ if there exists $f \in \mathscr S \prth{\R^2,\C}$ such that, for all $x \in \R^2, \vec f (x) = f(x) \vec V$.

\paragraph{Geometrical setup.}
We first introduce the geometrical setup.

\begin{assumption}[Transversality condition] \label{ass:transversality} The function $\m$ satisfies
\[
    \inf_{x \in \E} { \Big| \nabla \m (x) \Big| } > 0.
\]
\end{assumption}

Let $x_0 \in \E$.
In view of Assumption~\ref{ass:transversality}, we consider an arc-length parametrization $\mathbf{t}$ of the curve $\E$, such that $\mathbf{t}(0) = x_0$, so
\[
    \mathbf{t}'(s) \coloneqq \frac{{\nabla \m(\mathbf{t}(s))}^\perp}{\left| \nabla \m(\mathbf{t}(s)) \right|}, \quad s \in \R,
\]
where $^\perp$ denotes the counterclockwise rotation by $\pi/2$.
We also introduce the vector-valued function $\mathbf{n}$ defined by
\[
    \mathbf{n}(s) \coloneqq \frac{\nabla \m(\mathbf{t}(s))}{\left| \nabla \m(\mathbf{t}(s)) \right|}, \quad s \in \R.
\]
The curvature of the curve $\E$ at $\mathbf{t}(s)$ is the real number $\kappa(s)$, defined by Frenet's formula
\[
    \kappa(s) \coloneqq -\mathbf{t}''(s) \cdot \mathbf{n}(s).
\]
For all $s \in \R$, we set $r(s) \coloneqq \left| \nabla \m(\mathbf{t}(s)) \right|^{1/2}$ and $\displaystyle{\theta(s) \coloneqq \theta_0 + \int_0^s \kappa(\eta)d\eta}$, where $\theta_0 \in [0,2\pi)$ is chosen such that  $\displaystyle{\mathbf{n}(0) = \begin{pmatrix} -\sin \theta_0 \\ \cos \theta_0 \end{pmatrix}}$.
Then
\[
    \frac{\nabla \m(\mathbf{t}(s))}{\left| \nabla \m(\mathbf{t}(s)) \right|} = \begin{pmatrix} -\sin \theta(s) \\ \cos \theta(s) \end{pmatrix}.
\]

\paragraph{Generated edge states (Theorem 2 in~\cite{Edge-States} and Theorem A.1 in~\cite{Drouot}).}
Let $\vec f \in \mathscr S \prth{\R^2,\C^2}$ be independent of $\ep$.
We consider the following initial condition, where for all $\ep \in (0,1]$,
\begin{equation} \label{eq:initialconditionBBDFLW}
    \psi_0^\ep(x) \coloneqq \frac{1}{\sqrt\ep} \vec f\prth{ \frac{x-x_0}{\sqrt\ep} } = \mathrm{WP}^\ep_{\mathbf{t}(0),0}\left[ \vec{f} \right](x), \quad x \in \R^2.
\end{equation}
This initial wave packet is concentrated on the curve $\mathcal{C} \coloneqq \E \times \{(0,0)\} \subset \R^2_x \times \R^2_\xi$, which physically means that the wave packets carry no momentum.
Classically, this implies that they are not expected to move.

\bigskip

However, if $\displaystyle{ \prth{\psi_t^\ep}_{\ep > 0} }$ solves~\eqref{eq:Dirac} with the initial condition~\eqref{eq:initialconditionBBDFLW}, then there exists $T > 0$, such that, for all $t \in (0,T)$, the following holds uniformly for $\ep \in (0,1]$,
\begin{equation} \label{eq:solutionDrouot}
    \psi_t^\ep = \mathrm{WP}^\ep_{\mathbf{t}(t),0} \left[ F(t,\cdot) \vec V_{\theta(t)} \right] + \gO{\ep^{-1/4}}{\mathrm L^\infty\prth{\R^2,\C^2}} + \gO{\ep^{1/2}}{\L\prth{\R^2,\C^2}}.
\end{equation}
where $F \in \mathscr{C}^\infty \left( \R_t, \mathscr S \prth{\R^2,\C} \right)$ is determined by the initial data and $\displaystyle{\vec V_\theta \coloneqq \matp{e^{\frac{-\i\theta}{2}} \\ -e^{\frac{\i\theta}{2}} }}$.
Moreover,
\begin{enumerate}[label=(\roman*)]
    \item If $\prth{\psi_0^\ep}_{\ep>0}$ is oriented along $\vec V_{\theta_0}$, then the equality~\eqref{eq:solutionDrouot} holds for all $t>0$ without the $\gO{\ep^{-1/4}}{\mathrm L^\infty}$ term.
    \item If there exists $f \in \mathscr S \prth{\R,\C}$ such that, for all $\ep \in (0,1]$,
    \[
        \psi_0^\ep(x) \coloneqq \frac{r(0)^{1/2}}{\sqrt\ep} f \left( \frac{\left( \mathrm R_{\theta_0} \prth{x-x_0} \right)_1}{\sqrt\ep} \right) \exp \left( \frac{-r(0)^2 \left( \mathrm R_{\theta_0} \prth{x-x_0} \right)_2^2}{2\ep} \right) \vec V_{\theta_0}, \quad x \in \R^2,
    \]
    where $\mathrm R_\theta$ denotes the rotation of angle $\theta$ in $\R^2$, then, uniformly for $\ep \in (0,1]$ and for all $(t,x) \in \R \times \R^2$,
    \[
        \psi_t^\ep(x) = \frac{r_t^{1/2}}{\sqrt\ep} f \left( \frac{\left( \mathrm R_{\theta_t} \prth{x-\mathbf{t}(t)} \right)_1}{\sqrt\ep} \right) \exp \left( \frac{-r_t^2 \left( \mathrm R_{\theta_t} \prth{x-\mathbf{t}(t)} \right)_2^2}{2\ep} \right) \vec V_{\theta_t} \ + \ \gO{\ep^{1/2} \langle t \rangle}{}.
    \]
        Here, the variables $\left( \mathrm R_{\theta_t} \prth{x-\mathbf{t}(t)} \right)_1$ and $\left( \mathrm R_{\theta_t} \prth{x-\mathbf{t}(t)} \right)_2$ represent tangent and normal coordinates to the curve $\E$, respectively.
        They play different roles in the wave packet.
        In particular, for all $t \in [0,T)$, the wave packet part exhibits Gaussian behavior in the normal coordinate and remains localized near the point $(\mathbf{t}(t),(0,0))$ in the phase space.
    \item If $\prth{\psi_0^\ep}_{\ep>0}$ is oriented along $\vec V_{\theta_0}^\perp$, then $F = 0$.
        This implies that if the orientation is not aligned with $\vec V_{\theta_0}$, or if the decay in the normal variable is orthogonal in $\L$ to a Gaussian profile, then the $\mathrm L^\infty$ remainder is not necessarily negligible in $\L$-norm.
        Furthermore, the terms in the $\mathrm L^\infty$ remainder are not wave packets.
\end{enumerate}

\begin{SCfigure}[0.5][h]
    \caption{Example of geometric setup.}
    \includegraphics[width=0.6\textwidth]{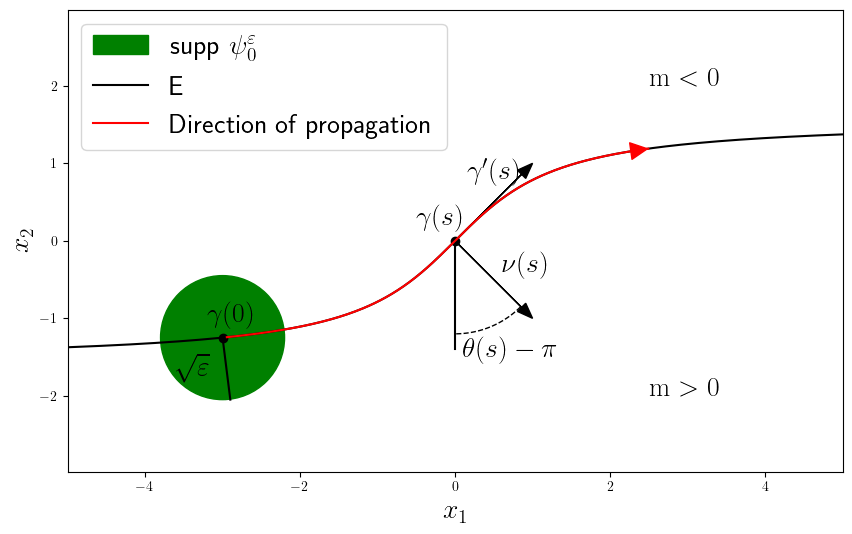}
\end{SCfigure}

In~\cite{Drouot}, the $\mathrm L^\infty$ remainder is described using a WKB expansion, for initial data that are wave packets, where the Hamiltonian is reduced to a normal form via a Fourier Integral Operator.
In~\cite{Bal}, a parametrix is constructed for wave packets used as initial conditions.

\bigskip

Our objective is twofold: we aim to describe the evolution of dispersive terms, and we also consider initial data in $\L\prth{\R^2,\C^2}$, not necessarily coherent states as in~\eqref{eq:initialconditionBBDFLW}.
We will focus on computing the Wigner measure of the solution to~\eqref{eq:Dirac} for initial data in $\L \left( \R^2,\C^2\right)$.

\subsection{Wigner transform and semiclassical measure} \label{subsec:wigner}

\paragraph{Wigner transform.}
In the present paper, instead of computing the asymptotics of the wave function itself, as in~\eqref{eq:solutionDrouot}, we focus on probability densities, in particular the position density defined by $\left|\psi^\ep_t(x)\right|_{\C^2}^2\diff x$.
Its physical interpretation is the probability of finding the particle at time $t$ in position $x$.
This approach generalizes the previous strategy and captures information on all quantities that are not negligible in $\L$ norm.

\bigskip

Since our problem requires an analysis in phase space, we go one step further and consider the {\it Wigner transform} of $\prth{\psi_t^\ep}_{\ep > 0}$ and its weak limits in the space of distributions, known as {\it Wigner measures}.
These were introduced by Wigner in~\cite{Wigner2} and~\cite{Wigner}.
The Wigner transform of $\displaystyle{(\psi_t^\ep)_{\ep > 0}}$ is defined as follows,
\[
    \mathrm W^\ep\crch{\psi_t^\ep}(x,\xi) \coloneqq \int_{\R^2} e^{\i\xi \cdot y} \psi_t^\ep \prth{x + \frac{\ep y}{2}} \otimes \overline{\psi_t^\ep \prth{x - \frac{\ep y}{2}}} \frac{\diff y}{(2\pi)^2}, \quad (x,\xi) \in \R^2 \times \R^2,
\]
where $v\otimes w$ denotes the $2 \times 2$ matrix $v\,^tw$ for $(v,w) \in \C^2 \times \C^2$.
The Wigner transform $\mathrm W^\ep\crch{\psi_t^\ep}$ belongs to $\L\prth{\R^4,\C^{2,2}} \cap \mathscr C ^0 \prth{\R^4,\C^{2,2}}$ and takes values in the set of Hermitian matrices.
Moreover, the position density can be formally recovered through the relation
\[
    \int_{\R^2} {\rm tr} \mathrm W^\ep\crch{\psi_t^\ep}(x,\xi) \diff\xi = \abs{\psi_t^\ep(x)}_{\C^2}^2.
\]
It is also possible to recover the momentum density by integrating over position rather than momentum.

\paragraph{Distributional Properties.}
The function $\mathrm W^\ep\crch{\psi_t^\ep}$ is a distribution on the phase space $\R^2_x \times \R^2_\xi$.
In fact,
\begin{equation} \label{eq:Wignerdistribution}
    \mathrm I^\ep_{\psi_t^\ep}(a) \coloneqq \int_{\R^{4}}{\rm tr} \Big( a(x,\xi) \mathrm W^\ep\crch{\psi_t^\ep}(x,\xi)\Big) \diff x\diff\xi = \Big\langle \Op(a)\psi_t^\ep , \psi_t^\ep \Big\rangle_{\L\prth{\R^2,\C^2}}, \quad a \in \mathscr C ^\infty_c (\R^4,\C^{2,2})
\end{equation}
where $\Op(a)$ denotes the semiclassical pseudodifferential operator with symbol $a$, obtained via Weyl quantization.
That is,
\[
    \Op(a)\psi_t^\ep(x) \coloneqq \int_{\R^4} a \prth{\frac{x+y}{2}, \xi} \psi_t^\ep(y) e^{\i\frac{\xi \cdot(x-y)}{\ep}} \frac{\diff y \diff\xi}{(2\pi\ep)^2}, \quad x \in \R^2.
\]
The Calder\'on-Vaillancourt Theorem~\cite{CV} asserts that the family of operators $\prth{\Op(a)}_{\ep > 0}$ is uniformly bounded with respect to $\ep$ in $\mathcal L\prth{\L\prth{\R^2,\C^2}}$:
there exists constants $C > 0$ and $N \in \N$ such that, for all $a = (a_{i,j})_{1\leqslant i,j\leqslant 2} \in \mathscr C^\infty_c \prth{\R^4,\C^{2,2}}$ and all $f \in \L\prth{\R^2,\C^2}$,
\begin{equation} \label{eq:Calderon_Vaillancourt}
    \norme{\Op(a)f}_{\L\prth{\R^2,\C^2}} \leqslant C \,\max_{1\leqslant i,j\leqslant 2}\, \max_{\abs{\alpha} + \abs{\beta} \leqslant N} \norme{\partial^\alpha_x \partial^\beta_\xi a_{i,j}}_{\mathrm L^\infty\prth{\R^4}} \norme{f}_{\L\prth{\R^2,\C^2}}.
\end{equation}
However, the Wigner transform is not a positive distribution.
Nevertheless, its weak limits are positive measures.
This follows from the weak G\aa rding inequality.
Adapting~\cite[Proposition 2.13]{Clotildecourse} to the matrix-valued case, we have: for all non-negative $a \in \mathscr C^\infty_c \prth{\R^4,\C^{2,2}}$ and for all $\delta > 0$, there exists a constant $C_\delta > 0$ such that,
\[
    \scal{\Op (a) f}{f}_{\L\prth{\R^2,\C^2}} \geqslant -\prth{\delta + C_\delta \ep}\norme{f}^2_{\L\prth{\R^2,\C^2}}, \quad f \in \L\prth{\R^2,\C^2}.
\]

\paragraph{Time-averaged Wigner measure.}
To account for the time dependence of the family $\prth{\psi^\ep_t}_{\ep > 0}$, we consider \textit{time-averaged Wigner measures}.
More precisely, we test the Wigner transform $\mathrm W^\ep[\psi^\ep_t]$ against separated time-dependent observables: for all $\Xi \in \mathscr C_c^\infty\prth{\R,\C}$, for all $a \in \mathscr C^\infty_c \prth{\R^4,\C^{2,2}}$,
\[
    \mathcal I^{\ep}_{\psi_t^{\ep}}(\Xi,a) \coloneqq \int_\R \int_{\R^{4}} \Xi(t) {\rm tr}_{\C^2}\Big( a(x,\xi) \mathrm W^\ep[\psi^\ep_t](x,\xi)\Big) \diff x\diff\xi\diff t \; = \; \int_\R \Xi(t) \scal{\Op(a)\psi^\ep_t}{\psi^\ep_t}_{\L\prth{\R^2,\C^2}} \diff t.
\]

\paragraph{Semiclassical measure.}
Our goal is to describe the solution $\prth{\psi_t^\ep}_{\ep > 0}$ to equation~\eqref{eq:Dirac} for initial data in $\L\prth{\R^2,\C^2}$.
At the very least, we can characterize the limit as $\ep$ goes to $0$ by constructing \textit{semiclassical measures}.
More precisely, there exists a vanishing sequence of positive numbers $\prth{\ep_k}_{k \in \N}$ and a positive matrix of Radon measures $\mu$ on $\R_t \times \R^2_x \times \R^2_\xi$, such that,
\begin{equation} \label{eq:semiclassical measure}
    \mathcal I^{\ep_k}_{\psi_t^{\ep_k}}(\Xi,a)  \; \underset{k \to +\infty}{\longrightarrow} \; \int_{\R^5} \Xi(t) {\rm tr}\Big( a(x,\xi) \diff\mu(t, x, \xi) \Big), \quad (a,\Xi) \in \mathscr C^\infty_c \prth{\R^4, \C^{2,2}} \times \mathscr C_c^\infty(\R,\C).
\end{equation}
The measure $\mu$ is called the \textit{time-averaged Wigner measure} of the family $\prth{\psi_t^{\ep_k}}_{k \in \N}$.
A proof of existence is recalled in Appendix~\ref{app:separable} for families taking values in a separable Hilbert space, as this will be a main topic of Section~\ref{sec:Two-scale semiclassical measures} (while observables and quantization will be precised in Section~\ref{sec:Two-scale Wigner observables}).
Moreover, thanks to the control
\[
    \abs{  \int_{\R^5} \Xi(t)a(x,\xi) \diff\mu(t, x, \xi) } \leqslant C \ \norme{\Xi}_{\mathrm L^1(\R)} \norme{\Op(a)}_{\mathcal L \prth{\L\prth{\R^2,\C^2}}},
\]
we can write $\diff\mu(t,x,\xi) = \diff\mu_t(x,\xi)\diff t$.
In other words, the measure $\mu$ is absolutely continuous in time with respect to the Lebesgue measure.
Our goal is to describe the propagation of these measures.

\subsection{Bulk propagation of the semiclassical measure}

It is known from \cite{GMMP, Clotildecourse} that any time-averaged Wigner measure of the family $\prth{\psi^\ep_t}_{\ep > 0}$, outside the curve $\mathcal C$, is completely determined by the corresponding measures of the initial data.
Let us first introduce some notations.
Let us consider the symbol of the Hamiltonian in the Weyl quantization defined as
\[
    \mathrm H(x,\xi) \coloneqq \begin{pmatrix}  \m(x) & \xi_1 - \i\xi_2 \\ \xi_1 + \i\xi_2 & - \m(x) \end{pmatrix}, \quad (x,\xi) \in \R^4.
\]
Notably, $\displaystyle{ \det \mathrm H(x,\xi) = - \m(x)^2 - |\xi|^2 }$ so that $\mathcal C = \{ (x,\xi) \in \R^4 \mid \det \mathrm H(x,\xi) = 0\}$.
Let $(x,\xi) \in \R^4 \setminus \mathcal C$.
Then, the matrix $\mathrm H(x,\xi)$ has two distinct eigenvalues, denoted by $\lambda_+(x,\xi)$ and $\lambda_-(x,\xi)$ defined by
\[
    \lambda_+(x,\xi) \ = \ -\lambda_-(x,\xi) \ \coloneqq \ \sqrt{ \m(x)^2 + |\xi|^2} \ > \ 0,
\]
and the projectors onto the corresponding eigenspaces are given by
\[
    \Pi_+(x,\xi) \coloneqq \mathrm{Id} - \frac{1}{\lambda_+(x,\xi)} \mathrm H(x,\xi), \quad \Pi_-(x,\xi) \coloneqq \mathrm{Id} - \frac{1}{\lambda_-(x,\xi)} \mathrm H(x,\xi).
\]
We recall the definition of the \textit{Poisson bracket}
\begin{equation} \label{eq:poissonbracket}
   \ens{f,g} \coloneqq \nabla_\xi f \nabla_x g - \nabla_x f \nabla_\xi g \; , \; (f,g) \in \mathscr C^1\prth{\R^2_x \times \R^2_\xi}^2.
\end{equation}
Following~\cite[section 6]{GMMP}, any semiclassical measure $\mu$ can be decomposed outside $\mathcal C$.

\begin{lemma}[Evolution of semiclassical measure outside $\mathcal C$] \label{lem:GMMP}
Let $\mu_t$ be a semiclassical measure of the family $\prth{\psi^\ep_t}_{\ep > 0}$ associated with the sequence $(\ep_n)_{n\in\N}$.
Then,
\[
    \diff\mu_t (x,\xi) = \Pi_+(x,\xi)\diff\mu_t^+ (x,\xi)\diff t + \Pi_-(x,\xi)\diff\mu_t^- (x,\xi)\diff t, \quad (x,\xi) \in \R^4 \setminus \mathcal C,
\]
where the scalar measures $\mu_t^\pm$ satisfy
\[
    \partial_t \mu_t^\pm = \ens{\mu_t^\pm , \lambda_\pm}, \quad \mu_0^\pm = \mathrm{tr}\Big( \Pi_\pm \mu_0 \Big), \\
\]
and $\mu_0$ is the semiclassical measure of the family $\prth{\psi^\ep_0}_{\ep > 0}$ associated with the sequence $(\ep_n)_{n\in\N}$.
\end{lemma}

This theorem implies that, for all initial data $\prth{\psi_0^\ep}_{\ep > 0}$ microlocalized around $(x_0,\xi_0) \notin \mathcal C$, the semiclassical measures evolve along the Hamiltonian trajectories associated with the functions $\lambda_+$ and $\lambda_-$.
These trajectories lie on the hypersurface $\det \mathrm H(x,\xi) = \det \mathrm H(x_0,\xi_0)$.
In other words, if the initial data are localized far from the curve $\mathcal C$, then the semiclassical measure never approaches the curve $\mathcal C$ in phase space.

\bigskip

Our aim in this paper is to describe the evolution of the measure $\mu$ above the curve $\mathcal C$.
As previously noted in~\cite{Edge-States} and~\cite{Drouot}, the concentration of $\prth{\psi_0^\ep}_{\ep > 0}$ is a two-scale problem involving the scale $\sqrt{\ep}$.
Therefore, we will adopt a two-scale Wigner approach at scale $\sqrt{\ep}$, as initiated in~\cite{FG2002} (see also~\cite{Clotildecourse, Nier, Miller}).
To achieve this objective, we work in a neighborhood of $\E$.

\subsection{Two-scale analysis above the interface}

\paragraph{Geometrical setup (2).}
To simplify the presentation, we will impose an assumption on $\m$ to provide a global tubular neighborhood of $\E$.
Let us denote
\[
    \mathrm I \coloneqq \left( -\frac{1}{2|| \kappa ||_{\infty}} , \frac{1}{2||\kappa ||_{\infty}}\right).
\]
We assume $\E$ is non-compact.
The following assumption ensures the existence of a global tubular neighborhood of size $|\mathrm{I}|$.

\begin{assumption} \label{ass:tubularneighborhood} The map
\begin{equation} \label{eq:phi}
    \begin{matrix} \Phi : & \R \times \mathrm I & \longrightarrow & \R^2 \\ & (s,y) & \longmapsto & \mathbf{t}(s) + y \mathbf{n}(s) \end{matrix}
\end{equation}
is a global diffeomorphism from $\R \times \mathrm I$ to
\[
    \Omega \coloneqq \left\{ \mathbf{t}(s) + y \mathbf{n}(s) | \ (s,y) \in \R \times \mathrm I \right\}.
\]
\end{assumption}

The set $\Omega$ is a global tubular neighborhood of the curve $\E$~\cite{Hirsch}.
Actually, this assumption is always satisfied locally.

\bigskip

Let us assume Assumption~\ref{ass:transversality} and Assumption~\ref{ass:tubularneighborhood}.
We will work in the normal geodesic coordinates, denoted $(s,y)$, defined by~\eqref{eq:phi}, and denote $(\sigma,\eta)$ as the dual variables of $(s,y)$, which play the role of $\xi$.
In these new coordinates, the curve $\mathcal C$ is reduced to $\ens{y = \sigma = \eta = 0}$.
Therefore, the measure $\mu$ over $\mathcal C$ can be rewritten as
\begin{equation} \label{eq:measure_change_variable}
    \diff \mu \diff t \mathds{1}_\mathcal C = \Phi^* \Big( \mathrm U_{\theta(s)}^\dagger \diff \rho_t \mathrm U_{\theta(s)} \Big) \diff t
\end{equation}
where $\rho_t$ is a semiclassical measure on $\R_s$ and the unitary operator $\mathrm U_\theta$ is the multiplication by a matrix defined by
\begin{equation} \label{eq:U_theta}
    \mathrm U_\theta \coloneqq \frac{1}{\sqrt 2} \begin{pmatrix} e^{ \i\frac{\theta}{2}} & -e^{- \i\frac{\theta}{2}} \\ e^{ \i\frac{\theta}{2}} & e^{- \i\frac{\theta}{2}} \end{pmatrix}, \quad \theta \in \R.
\end{equation}

\paragraph{Two-scale quantization.}
Based on the geometric assumptions on $\mathcal C$, we introduce a two-scale quantization procedure.
We follow a two-scale Wigner approach at scale $\sqrt{\ep}$, as initiated in~\cite{FG2002} (see also~\cite{Clotildecourse}, \cite{Nier} and~\cite{Miller}).
The approach using tubular coordinates is inspired by~\cite{Fab}.

\bigskip

We say that $a \in \mathscr C^\infty \Big( \R^4 \times \R^3, \C^{2,2} \Big)$ is a two-scale observable if, and only if,
\begin{enumerate}[label=(\roman*)]
    \item there exists a compact set $K \subset \R^4$ such that for all $z \in \R^3$, $a(\cdot,\cdot,z)$ is supported in $K$,
    \item there exists a constant $R_0 > 0$ and $a_\infty \in \mathscr C^\infty \Big( \R^4 \times \S^2, \C^{2,2} \Big)$, such that for all $|z| > R_0$, $a(\cdot,\cdot,z) = a_\infty\left(\cdot,\cdot,\frac{z}{|z|}\right)$.
\end{enumerate}
We denote by $\mathcal{A}$ the set of two-scale observables.
We define a second-scale quantization near the interface $\mathcal C$ for these observables by
\begin{equation} \label{eq:secondquantization}
\begin{array}{ccccc}
    \Opp & : & \mathcal{A} & \longrightarrow & \mathcal L\big( \L\prth{\R^2,\C^2}\big) \\
     & & a & \longmapsto & \mathrm{Op}_1 \left( a \left( s,y , \ep \sigma, \ep \eta , \frac{y}{\sqrt\ep}, \sqrt\ep \sigma, \sqrt\ep \eta \right) \right)
\end{array}
\end{equation}
The quantization map $\Opp$ associates to $a$ a bounded operator on $\L\prth{\R^2,\C^2}$ according to the Calder\'on-Vaillancourt Theorem~\cite{CV} as we will see in Section~\ref{sec:Two-scale Wigner observables}.
We can identify our observables in $\mathscr C^\infty_c \prth{\R^4, \C^{2,2}}$ as a subset of $\mathcal{A}$.
With an abuse of notation, for all $a \in \mathscr C^\infty_c \prth{\R^4, \C^{2,2}}$, $a$ belongs to $\mathcal{A}$ with $a_\infty = a$ and moreover
\[
    \Opp(a) = \Op(a).
\]
By analyzing the limits~\eqref{eq:semiclassical measure} for observables supported on $\Omega$, we can replace $\Op(a)$ by $\Opp(a)$, which will gives us (according to Section~\ref{sec:Two_scale_analysis}) the following decomposition where $\S^2$ denotes the unit sphere of $\R^3$, for all $(s,\sigma) \in \R^2$, $a^W (s;0; y,\sigma,\D_y)$ denotes the bounded operator on $\L\prth{\R_y, \C^2}$ obtained by the Weyl quantization of $(y,\eta)\mapsto a(s;0;y,\sigma,\eta)$ and $\tr$ denotes the trace in $\mathcal L \prth{ \L\prth{\R,\C^2} }$.

\begin{theorem}[Time-averaged two-scale semiclassical measure above the interface] \label{thm:shape_two_scale_measure}
Let $\displaystyle{(\psi^\ep_t)_{\ep > 0}}$ be a solution to \eqref{eq:Dirac} with normalized initial condition in $\L \prth{\R^2,\C^2}$.
Then, there exist a vanishing sequence of positive numbers $\displaystyle{(\ep_k)_{k \in \N}}$ and
\begin{enumerate}[label=(\roman*)]
    \item Two measurable maps $t\mapsto \diff\nu^t_\infty$ and $t\mapsto \diff\nu^t$, valued in the set of non-negative scalar measures on $\mathcal C \times \S^2$ and $\R^2$, respectively,
    \item A $\diff\nu^t_\infty\otimes \diff t $-measurable map 
    $(t,s,\omega)\mapsto M^t_\infty(s,\omega) \in \C^{2,2}$, valued in the set of Hermitian trace $1$ positive matrices,
    \item A $\diff\nu^t(s,\sigma)\otimes \diff t $ measurable map $(t,s,\sigma)\mapsto M^t(s,\sigma)$, valued in trace-one positive operators on $\L\left(\R, \C^2\right)$,
\end{enumerate}
such that the semiclassical measure $\rho_t$ defined in~\eqref{eq:measure_change_variable} associated with $\displaystyle{(\ep_k)_{k \in \N}}$ satisfy, for all $s \in \R$,
\[
    \rho_t(s) = \int_{\S^2} \mathrm{tr} \Big( a_\infty (s;0; \omega) M_\infty^t(s,\omega) \Big) \diff\nu_\infty^t(s,\omega) + \int_{\R} \tr \Big( a^W (s;0;y,\sigma,\D_y) M^t(s,\sigma) \Big) \diff\nu^t(s,\sigma).
\]
\end{theorem}

This Theorem is a direct consequence of two results.
The first one is the existence of two-scale semiclassical measures, which is Theorem~\ref{thm:existence_two_scale_measure} and the second one is the link between semiclassical measures above the interface, which is Lemma~\ref{lem: normal form measure}.

\begin{rem*}[Geometrical point of view] In the parametrization $(s,y)$ of $\R\times \mathrm I$, the variables $(s,y,\sigma,\eta)$ parametrize $\mathrm T^*\R^2$.
Moreover, the variable $s$ parametrizes the curve $\E \subset \R^2$.
The variable $y$ parametrizes the normal fiber to the curve $\E$ at the point $\mathbf{t}(s)$.
In the computation of the trace $\tr$, the space $\L(\R_y,\C^2)$ can be understood as the space $\L(N_{\mathbf{t}(s)}\E,\C^2)$, where $N_{\mathbf{t}(s)} \E$ denotes the normal bundle of the curve $\E$ at the point $\mathbf{t}(s)$.
The variable $\sigma$ denotes the variable of the fiber of the cotangent bundle $\mathrm T^*\E$ over $\E$.
Specifically, for $\mathbf{t}(s) \in \E$ and $\mathrm T_s \E$, the tangent to $\E$ at $s$, the fiber of the cotangent bundle $\mathrm T^*\E$ over $\mathbf{t}(s)$ is $(\mathrm T_s\E)^*\sim\R$, which is one-dimensional space.
The measure $\nu^t$ defines a measure on $\mathrm T^*\E$, the cotangent bundle of the curve $\mathcal{E}$.

\medskip

Since  $\mathcal C=\{y=\eta=\sigma=0\}\subset \mathrm T^*\R^2$, the set $\mathcal C\times \S^2$
coincides with the spherical normal bundle of $\mathcal C$ denoted by $\mathrm S\mathcal C$.
The fiber of $\mathrm S \mathcal C$ over a point $\rho\in\mathcal C$ is defined as follows.
Let $\mathrm T_\rho\mathcal C$ be the tangent space to $\mathcal C$ at $\rho$, and define $\mathrm N_\rho \mathcal C$ at the subspace of $\prth{ \mathrm T_\rho \prth{ \mathrm T^*\R^2}}^*$ consisting of linear forms on $\mathrm T_\rho \prth{\mathrm T^*\R^2}$ that vanish on $\mathrm T_\rho \mathcal C$.
Then $\mathrm S_\rho\mathcal C$ is obtained by taking the quotient of $\mathrm N_\rho\mathcal C$ by the action of $\R_{\geqslant 0}$ by homotheties.
\end{rem*}

\begin{def*}[Two-scale semiclassical measure] 
We denote by $M\diff\nu$ the pair $\prth{M, \nu}$, where $M$ is a positive trace-one operator that is $\diff\nu$-measurable, and $\nu$ is a non-negative measure.
We denote by $M_\infty\diff\nu_\infty$ the pair $\prth{M_\infty, \nu_\infty}$ with $M_\infty$ a Hermitian trace one positive matrix $\diff\nu_\infty$-measurable and $\nu_\infty$ a non-negative measure.

\medskip

We call two-scale semiclassical measures, associated with the concentration of the family $\displaystyle{\prth{\psi_t^\ep}_{\ep > 0}}$ near the curve $\mathcal C$ at scales between $\sqrt\ep$ and $\delta$, and at scale $\sqrt\ep$, respectively, the pair $\displaystyle{\prth{M_\infty\diff\nu_\infty, M\diff\nu}}$.
We say that $M\diff\nu$ is a two-scale semiclassical measure at finite distance, while $M_\infty\diff\nu_\infty$ is a two-scale semiclassical measure at infinity.
\end{def*}

\begin{SCfigure}[0.5][h]
    \caption{Illustration of two-scale separation of space before passing to limit.}
    \includegraphics[width=0.6\textwidth]{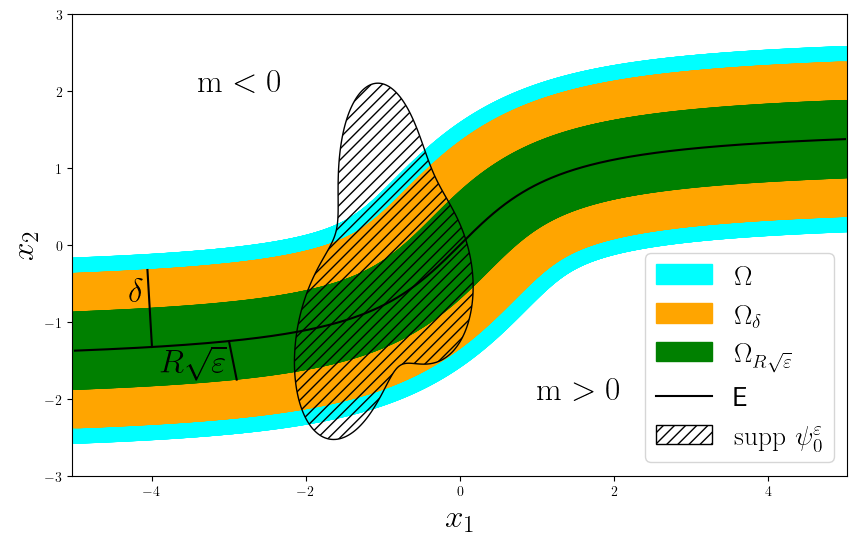}
\end{SCfigure}

\subsection{Main result : edge propagation of the semiclassical measure} \label{sec:main_result}

Let us specify the effective Hamiltonians associated with the different concentration regimes.
As in Lemma~\ref{lem:GMMP}, the eigenvalues and eigenprojectors of these Hamiltonians will play a crucial role in the propagation.

\paragraph{Concentration at finite distance.}
Let us consider the symbol of the principal term of the Hamiltonian, depending only on the tangential variable $s$, while retaining its operator character with respect to the normal variable $y$; namely, the operator acting on $\L\prth{\R_s,\L\prth{\R_y,\C^2}}$ of symbol
\begin{equation} \label{eq:Hamiltonian_pp_scale_ep}
    \mathrm T_\E(s,\sigma) \coloneqq \matp{ \sigma & y|\nabla \m (\mathbf{t}(s))| - \i\D_y \\ y|\nabla \m (\mathbf{t}(s))| + \i\D_y & -\sigma}, \quad (s,\sigma) \in \R^2.
\end{equation}
A complete spectral analysis of this operator is provided in Appendix~\ref{app:spectral_theory}.
Its eigenvectors depend on the Hermite functions $\prth{\mathfrak{h}_m}_{m \in \N_{\geqslant 0}}$, which form an Hilbertian basis of $\L\prth{\R,\C}$ and satisfy, for all $y \in \R$ and all $m \in \N_{\geqslant 0}$,
\begin{equation} \label{eq:Hermite_function}
    -\frac{\diff^2 \mathfrak{h}_m}{\diff y^2} (y) + y^2 \mathfrak{h}_m (y) = (2m + 1)\mathfrak{h}_m(y).
\end{equation}
Let us consider, for all $(s,y,\sigma) \in \R^3$,
\begin{align}
    g_0^s (y) \coloneqq \frac{1}{\sqrt{r(s)}}\begin{pmatrix} \mathfrak{h}_0\big(r(s)y\big) \\ 0 \end{pmatrix}, \; \; & \; \; g_n^{s,\sigma} (y) \coloneqq \frac{1}{\sqrt{2r(s)}} \sqrt{ 1 - \frac{\sigma}{\lambda_n(s,\sigma)}} \begin{pmatrix} \frac{\sqrt{2|n|} r(s)}{\lambda_n(s, \sigma) - \sigma} \mathfrak{h}_{|n|} \big(r(s)y\big) \\ \mathfrak{h}_{|n|-1}\big(r(s)y\big)\end{pmatrix}, \label{eq:eigenvector_scale_ep} \\
    \Pi_0 (s) \coloneqq g_0^s \otimes g_0^s, \; \; & \; \; \Pi_n (s,\sigma) \coloneqq g_n^{s,\sigma} \otimes g_n^{s,\sigma}, \label{eq:proj_scale_ep}
\end{align}
where
\begin{equation} \label{eq:lambda_scale_ep}
    \lambda_0(\sigma) \coloneqq \sigma, \; \; \; \; \lambda_n(s,\sigma) \coloneqq \sgn(n) \sqrt{\sigma^2 + 2|n||\nabla \m (\mathbf{t}(s))|}.
\end{equation}
Although the functions $g_0^s, \lambda_0$ and $\Pi_0$ do not depend on $\sigma$, we will also denote them by $g_0^{s,\sigma}, \lambda_0(s,\sigma)$ and $\Pi_0(s,\sigma)$.

\medskip

Since we retain the operator character in the variable $y$, we are naturally led to consider semiclassical measure on the separable Hilbert space $\mathcal H = \L\left(\R_y, \C^2\right)$.

\paragraph{Concentration at infinity.}
We consider the symbol of the Hamiltonian corresponding to the second-scale quantization~\eqref{eq:secondquantization} at infinity, defined by
\begin{equation} \label{eq:Hamiltonian_pp_scale_delta}
    \mathrm T^\infty_\E (s, z) \coloneqq \matp{ z_\sigma & z_y r(s)^2 - \i z_\eta \\ z_y r(s)^2 + \i z_\eta & -z_\sigma}, \quad (s,y) \in \R^2, \ z \coloneqq (z_y, z_\sigma, z_\eta) \in \R^3,
\end{equation}
which has two eigenprojectors
\begin{equation} \label{eq:proj_infty}
    \Pi_\pm^\infty (s,z) \coloneqq \frac{1}{2} \mathrm{Id} \pm \frac{1}{2 \sqrt{ z_\sigma^2 + z_y^2r(s)^4 + z_\eta^2 }} \matp{ z_\sigma & z_y r(s)^2 - \i z_\eta \\ z_y r(s)^2 + \i z_\eta & -z_\sigma}.
\end{equation}

\paragraph{Main result.}
We now describe the semiclassical measures of Theorem~\ref{thm:shape_two_scale_measure} for solutions of equation~\eqref{eq:Dirac}.
Their semiclassical measures and their evolution outside $\mathcal C$ are already well understood.
Therefore, we focus on the two-scale semiclassical measures.
By applying a method similar to that of~\cite{GMMP} to the operators $\mathrm T_\E$ and $\mathrm T^\infty_\E$ in normal geodesic coordinates, together with a two-scale approach, we obtain the following main result.

\begin{theorem}[Propagation of the two-scale semiclassical measures] \label{thm:main_result} 
Assume Assumption~\ref{ass:transversality} and Assumption~\ref{ass:tubularneighborhood}.
Let $\displaystyle{ \prth{M^t_\infty \diff\nu^t_\infty , M^t \diff\nu^t} }$ be a pair of two-scale semiclassical measures associated with the family $\prth{\psi^\ep_t}_{\ep > 0}$, solution to~\eqref{eq:Dirac} with normalized initial condition.
Then the following properties hold.
\begin{enumerate}
    \item There exists a family of non-negative scalar Radon measures $\prth{\nu_n^t}_{n \in \Z}$ such that
\begin{equation} \label{eq:measure_scale_ep}
M^t \diff \nu^t = \sum_{n \in \Z} \Pi_n \diff \nu_n^t.
\end{equation}
Moreover, for all $n\in\Z$, $\nu_n^t$ satisfies the following equation
\begin{equation} \label{eq:propagation_measure_scale_ep}
    \left\{
    \begin{array}{rll}
    \partial_t \nu_n^t & = & \ens{\nu_n^t, \lambda_n} \\
    \nu_n^t |_{t=0} & = & \nu_{n,0} 
    \end{array}
    \right. ,
\end{equation}
where $\nu_{n,0}= \mathrm{tr} \Big( \Pi_nM^0 \Big) \nu^0$, and $M^0 \diff\nu^0 $ is a two-scale semiclassical measure associated with the initial data.
    \item There exists a family of non-negative scalar Radon measures $\nu^t_\pm$ such that
\begin{equation} \label{eq:measure_scale_delta}
    M^t_\infty \diff\nu^t_\infty = \Pi_+^\infty \diff\nu^t_+ + \Pi_-^\infty \diff\nu^t_-.
\end{equation}
Moreover, $\nu_\pm^t$ satisfies the invariance equation $\mathrm{div}_\omega \prth{\vec{V}^\infty(s,\omega) \nu_\pm^t(s,\omega)} = 0$, with
\begin{equation} \label{eq:defV}
    \vec{V}^\infty (s,\omega) \coloneqq \matp{ \prth{\omega_y\omega_\eta \prth{1-r(s)^4} - 1} \omega_\eta \\ \prth{1-r(s)^4} \omega_y\omega_\eta \omega_\sigma \\ \prth{\omega_y\omega_\eta \prth{1-r(s)^4} +r(s)^4} \omega_y }, \quad s \in \R, \quad \omega \coloneqq \prth{\omega_y, \omega_\sigma, \omega_\eta} \in \S^2.
\end{equation}
\end{enumerate}
\end{theorem}

\begin{SCfigure}[0.6][h]
    \caption{Illustration of the propagation of semiclassical measures for general initial data.}
    \includegraphics[width=0.6\textwidth]{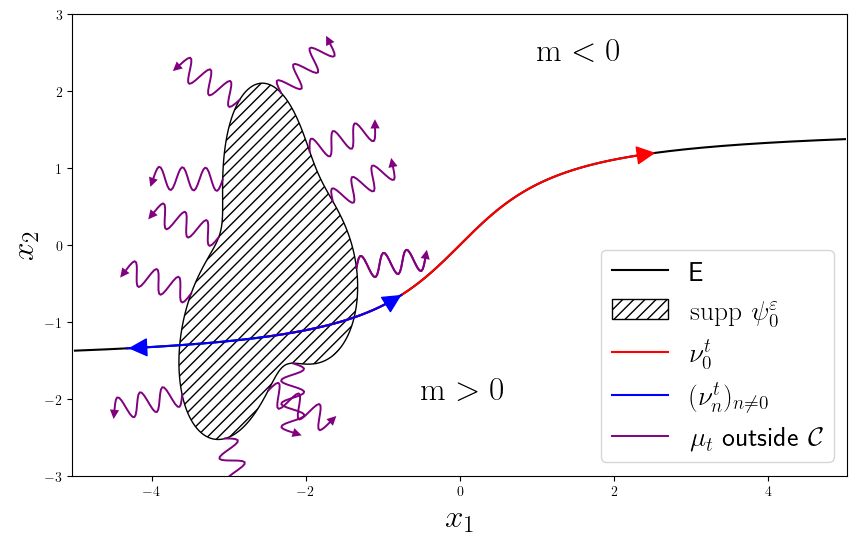}
\end{SCfigure}

\begin{rem*}[Continuity of the two-scale semiclassical measure at finite distance] From the evolution~\eqref{eq:propagation_measure_scale_ep}, it follows that the map $t \mapsto M^t\diff\nu^t$ is continuous and entirely determined by the evolution equation together with the initial data.
\end{rem*}

\begin{rem*}[Value at $t=0$ of the two-scale semiclassical measure at finite distance]
The value of $M^t\diff\nu^t$ at $t=0$ is not necessarily the two-scale semiclassical measure at finite distance of the initial condition.
As we will see in the proof, $M^t(s,\sigma)$ commutes with the operator $\mathrm T_\E(s,\sigma)$ due to the PDE satisfied by the solution near the interface.
Nevertheless, there exist initial conditions whose two-scale semiclassical measure at finite distance does not commute with $\mathrm T_\E(s,\sigma)$.
An example of this will be discussed in Remark~\ref{rem:difference_zero_measure}.
\end{rem*}

\begin{rem*}
Retaining the notation of the previous theorem, the measure $\nu_0^t$ corresponds to the wave packet associated with the solution studied in~\cite{Edge-States}.
Similarly, the measures $\prth{\nu_n^t}_{n \in \Z_{\neq 0}}$ correspond to the $\mathrm L^\infty$-remainder part of the solution in~\cite{Drouot}.
In particular, the mass of the edge mode is given by
\[
    \int_{\R^2} \diff\nu_{0,0}(s,\sigma).
\]
\end{rem*}

\subsection{Applications} \label{sec:applications}

Our main result describes the evolution of the semiclassical measure for arbitrary initial data, without any restriction on its structure.
A direct consequence is the following asymptotic description of the position density of the solution to equation~\eqref{eq:Dirac} for initial data given by wave packets concentrating on the curve $\mathcal C$.

\begin{corollary}[Evolution consequence] \label{cor:conservation}
Retaining the notations of Theorem~\ref{thm:main_result}, if
\[
    \int_{\R^4 \backslash \mathcal{C}} \mathrm{tr} \Big( \diff\mu_0(x,\xi) \Big) + \int_{\R^2} \tr \Big( M^0(s,\sigma) \Big) \diff\nu^0(s,\sigma) = 1,
\]
with $\mu_0$ a semiclassical measure of $\prth{\psi^\ep_0}_{\ep > 0}$ associated with the same vanishing sequence $\prth{\ep_n}_{n \in \N}$ as the two-scale semiclassical measure $\prth{M^0_\infty\diff\nu^0_\infty, M^0\diff\nu^0}$, it follows that, for almost all $t$,
\[
    M^t_\infty (s,\omega) \diff\nu^t_\infty (s,\omega) = 0.
\]
\end{corollary}

The evolution of the two-scale semiclassical measure in normal geodesic coordinates allows us to describe the general evolution of any normalized initial condition for equation~\eqref{eq:Dirac}.
This yields the following corollary, which is proved in Section~\ref{sec:proof_applications}.

\begin{corollary}[Propagation of semiclassical measure] \label{cor:goal_reached}
Assume Assumption~\ref{ass:transversality} and Assumption~\ref{ass:tubularneighborhood}.
Let $\prth{\psi^\ep_t}_{\ep > 0}$ be solutions to~\eqref{eq:Dirac} with normalized initial condition $\prth{\psi^\ep_0}_{\ep > 0}$ in $\L\prth{\R^2,\C^2}$.
Then, there exist a vanishing sequence of positive numbers $\displaystyle{(\ep_k)_{k \in \N}}$ and
\begin{itemize}
	\item[$\bullet$] a semiclassical measure $\mu_t$, evolving as described in Lemma~\ref{lem:GMMP},
	\item a family of non-negative scalar Radon measures $\prth{\nu_n^t}_{n \in \Z}$ satisfying~\eqref{eq:propagation_measure_scale_ep},
	\item[$\bullet$] a family $\nu^t_\pm$ of non-negative scalar Radon measures on $\mathcal C \times \S^2$ such that satisfying the invariance equation $\mathrm{div}_\omega \prth{\vec{V}^\infty(s,\omega) \nu_\pm^t(s,\omega)} = 0$, with $\vec{V}^\infty$ defined in~\eqref{eq:defV},
\end{itemize}
such that, for all $(a,\Xi) \in \mathscr C_c^\infty\prth{\R^4,\C^{2,2}} \times \mathscr C_c^\infty(\R,\C)$,
\[
    \mathcal I^{\ep_k}_\psi(\Xi,a) \underset{k \to + \infty}{\longrightarrow} \ \int_\R \int_{\R^4 \backslash \mathcal{C}} \Xi(t) \mathrm{tr} \Big( a (x,\xi) \diff\mu_t(x,\xi)\Big) \diff t + \int_\R \int_\R \Xi(t) \mathrm{tr} \Big( a \big( \mathbf{t}(s),0 \big) \diff\rho_t(s) \Big) \diff t,
\]
with, for all $s \in \R$,
\[
    \rho_t(s) = \sum_\pm \int_{\S^2} \mathrm U_{\theta(s)}^\dagger \Pi^\infty_\pm (s,\omega) \mathrm U_{\theta(s)} \diff\nu^t_\pm (s,\omega) + \frac{1}{2} \sum_{n \in \Z} \int_\R \matp{ 1 & \frac{-\sigma}{\lambda_n(s,\sigma)} e^{ -\i\theta(s)} \\ \frac{-\sigma}{\lambda_n(s,\sigma)}e^{ \i\theta(s)} & 1} \diff\nu^t_n(s,\diff\sigma).
\]
\end{corollary}

\paragraph{Link between~\cite{Drouot}, \cite{Edge-States} and Theorem~\ref{thm:main_result}}
A direct corollary of Theorem~\ref{thm:main_result} is the following.

\begin{corollary} \label{cor:linkDrouot}
Assume Assumption~\ref{ass:transversality} and Assumption~\ref{ass:tubularneighborhood}.
Let $\vec f \in \L\prth{\R^2,\C^2}$.
If $\displaystyle{ \prth{\psi_t^\ep}_{\ep > 0} }$ solves~\eqref{eq:Dirac} with, for all $\ep \in (0,1]$,
\[
    \psi_0^\ep(x) = \frac{1}{\sqrt\ep} \vec f \prth{ \frac{x-x_0}{\sqrt\ep} }, \quad x \in \R^2,
\]
then there exist a vanishing sequence of positive numbers $\displaystyle{(\ep_k)_{k \in \N}}$ and a family of measurable maps $\displaystyle{t\mapsto \prth{\nu^t_n}_{n \in \Z}}$ valued in the set of non-negative scalar Radon measures on $\R^2$, such that for all $(a,\Xi) \in \mathscr C_c^\infty \prth{ \R^4 ,\C^{2,2}} \times \mathscr C_c^\infty( \R ,\C)$,
\[
    \int_\R \Xi(t) \scal{\Op(a) \psi^{\ep_k}_t}{\psi^{\ep_k}_t}_{\L\left(\R^2,\C^2\right)} \diff t \underset{k \to + \infty}{\longrightarrow} \int_\R \int_\R \Xi(t) \mathrm{tr} \Big( a \big( \mathbf{t}(s),0 \big) \diff\rho_t(s) \Big) \diff t
\]
with,
\[
    \rho_t(s) = \frac{1}{2} \sum_{n \in \Z} \int_\R \matp{ 1 & \frac{-\sigma}{\lambda_n(s,\sigma)} e^{ -\i\theta(s)} \\ \frac{-\sigma}{\lambda_n(s,\sigma)}e^{ \i\theta(s)} & 1} \diff\nu^t_n(s,\diff\sigma), \quad s \in \R,
\]
where, for all $n\in\Z$, $\nu_n^t$ satisfies~\eqref{eq:propagation_measure_scale_ep} and $\nu_n^0$ is determined by the initial condition.
\end{corollary}

With the notation of Corollary~\ref{cor:linkDrouot}, the measure $\nu_0^t$ corresponds to the wave packet part of the solution in~\cite{Edge-States}, i.e., the edge mode generated by the initial data~\eqref{eq:initialconditionBBDFLW}.
Similarly, the measures $\prth{\nu_n^t}_{n \in \Z_{\neq 0}}$ correspond to the $\mathrm L^\infty$-remainder part of the solution in~\cite{Drouot}.
Notably, if $\nu_n^0 = 0$ for all $n \in \Z_{\neq 0}$ and $\nu_0^0 = 1$, then, for all $n \in \Z_{\neq 0}$,
\[
    \nu_0^t(s,\sigma) = \delta(s-t) \delta(\sigma), \quad \nu_n^t(s,\sigma) = 0, \quad (t,s,\sigma) \in \R^3,
\]
where $\delta$ denotes the Dirac mass at $0$.

\bigskip

Corollary~\ref{cor:linkDrouot} provides a description of the evolution of the position density for initial data given by wave packets.
However, Theorem~\ref{thm:main_result} allows us to consider arbitrary data in $\L\prth{\R^2,\C^2}$ concentrating near the curve $\mathcal C$, and thus to describe their asymptotic position-density behavior.

\begin{corollary} \label{cor:below_estimate}
Assume Assumption~\ref{ass:transversality} and Assumption~\ref{ass:tubularneighborhood}.
Let $\displaystyle{ \prth{\psi_0^\ep}_{\ep > 0} }$ be a uniformly bounded family in $\L\prth{\R^2,\C^2}$, let $\nu$ be a non-negative scalar Radon measure on $\R_t \times \R_x^2$, and let $\displaystyle{(\ep_k)_{k \in \N}}$ be a vanishing sequence such that, for all $(a,\Xi) \in \mathscr C_c^\infty \prth{ \Omega ,\C} \times \mathscr C_c^\infty( \R ,\C)$,
\[
    \int_\R \int_\Omega \Xi(t) a(x) \abs{\psi_t^{\ep_k}(x)}^2 \diff x \diff t \underset{k \to +\infty}{\longrightarrow} \int_\R \int_\Omega \Xi(t) a(x) \diff\nu(t,x),
\]
where $\displaystyle{ \prth{\psi_t^\ep}_{\ep > 0} }$ solves~\eqref{eq:Dirac} with $\displaystyle{ \prth{\psi_0^\ep}_{\ep > 0} }$ as initial condition.
Then,
\[
    \nu \mathds{1}_\E \geqslant \sum_{n \in \Z} \Phi^* \left( \int_\R \nu^t_n(\cdot,\diff\sigma) \right) \mathds{1}_\E \diff t
\]
where $\prth{\nu_n^t}_{n \in \Z}$ is defined by~\eqref{eq:propagation_measure_scale_ep}.
\end{corollary}

\noi This corollary provides a lower bound for the position density at time $t$ for the evolution of arbitrary initial data concentrating along the curve $\E$.

\subsection{Organization of the paper}

The proof of Theorem~\ref{thm:main_result} relies on two main arguments: the decomposition of the two-scale semiclassical measure above the curve $\mathcal C$ for arbitrary functions (Theorem~\ref{thm:shape_two_scale_measure}) and the decomposition of each resulting measure (Theorem~\ref{thm:main_result}) that takes into account equation~\eqref{eq:Dirac}.

\bigskip

First, to study the evolution of the measure $\mu_t$ near the curve $\mathcal C$, we perform a two-scale analysis in a neighborhood of $\mathcal C$.
Geometrically, $\mu_t$ is a measure on the cotangent space $\mathrm T^*\R^4$ and, as such, satisfies the associated geometric invariance properties.
In particular, this allows us to perform changes of coordinates.
We work locally in a tubular neighborhood of $\mathcal C$, where we can use the normal geodesic coordinates introduced earlier.
These coordinates $(s,y)$ are then used to transform the Hamiltonian, as detailed in Section~\ref{sec:normal geodesic coordinates}.
In Section~\ref{sec:Two-scale Wigner observables}, we separate and highlight two-scale Wigner observables associated with the different concentration regimes.
In Section~\ref{sec:Two-scale semiclassical measures}, we prove the existence of two-scale Wigner measures and describe them for arbitrary families.
A rescaling is then applied in Section~\ref{sec:rescaling} to take into account the two-scale Wigner approach of the Hamiltonian.
We will end up with a principal operator of the following form
\[
    \sqrt\ep \matp{ \sqrt\ep \D_1 & x_2 \widetilde{m}(x_1) + i \D_2 \\ x_2 \widetilde{m}(x_1) - i \D_2 & -\sqrt\ep \D_1}.
\]
This operator is $\sqrt\ep$-semiclassical with respect to $x_1$ only which justifies our two-scale approach.

\bigskip

Secondly, we describe each resulting measure and then prove our main result and its applications in Section~\ref{sec:proof}.
In fact, in Section~\ref{sub:Gamma}, we describe the evolution of the measure at finite distance $M^t\diff\nu^t$ and in Section~\ref{sub:M}, we describe the evolution of the measure at infinity $M_\infty^t\diff\nu_\infty^t$.
The decomposition of the measure $\mu_t$ as a sum of eigenprojectors associated with scalar measures is a consequence of the non-commutativity of matrices since
\[
    \int_\R \Xi(t) \i \sqrt\ep \frac{\mathrm d}{\mathrm{dt}} \scal{\Op(a)\psi^\ep_t}{\psi^\ep_t}_{\L\prth{\R^2,\C^2}} \diff t = \int_\R \Xi(t) \scal{ \Op \Big( \crch{a,\mathrm{H}} \Big) \psi^\ep_t}{\psi^\ep_t} \diff t + \mathcal{O}(\sqrt{\ep}).
\]
The evolution equation satisfied by the scalar measures then follows from the analysis of the next-order terms in the expansion in $\sqrt{\ep}$ which is
\[
    \int_\R \Xi(t) \i \sqrt\ep \frac{\mathrm d}{\mathrm{dt}} \scal{\Op(a)\psi^\ep_t}{\psi^\ep_t}_{\L\prth{\R^2,\C^2}} \diff t = \frac{\sqrt\ep}{2\i} \int_\R \Xi(t) \scal{ \Opd \Big( \ens{a,\mathrm{H}} - \ens{\mathrm{H}, a} \Big) \psi^\ep_t}{\psi^\ep_t} \diff t.
\]
We finally conclude by proving Theorem~\ref{thm:main_result} in Section~\ref{sec:proof_main} and by proving applications in Section~\ref{sec:proof_applications}.

\bigskip

In Appendix~\ref{app:Dirac}, we first apply the main properties of the Wigner transform (introduced in Section~\ref{subsec:wigner}) to the solution of our problem.
In Appendix~\ref{app:separable}, we recall the proof of existence of semiclassical measure for families valued in a separable Hilbert space.
This is necessarily since we study families valued in $\C$, $\C^2$ and $\L\prth{\R,\C^2}$.
In Appendix~\ref{app:spectral_theory}, we give a complete spectral analysis of~\eqref{eq:Hamiltonian_pp_scale_ep}.

\subsection{Notations}

We denote by $\displaystyle{ \prth{\sigma_j}_{1 \leqslant j \leqslant 3} }$ the Pauli matrices, defined as
\[
    \sigma_1 \coloneqq \begin{pmatrix} 0 & 1 \\ 1 & 0 \end{pmatrix}, \; \; \sigma_2 \coloneqq \begin{pmatrix} 0 & -\i \\ \i & 0 \end{pmatrix}, \; \; \sigma_3 \coloneqq \begin{pmatrix} 1 & 0 \\ 0 & -1 \end{pmatrix}.
\]
We denote $\H$ the operator
\[
    \H \coloneqq \begin{pmatrix}  \m(x) & \ep \D_1 - \i\ep \D_2 \\ \ep \D_1 + \i\ep \D_2 & - \m(x) \end{pmatrix} = \m (x)\sigma_3 + \ep \D_1 \sigma_1 + \ep \D_2 \sigma_2.
\]
Let $\mathcal H$ be a separable Hilbert space and $n$ an integer.
We denote the usual inner product by
\[
    \scal{f}{g}_{\L\prth{\R^n,\mathcal H}} \coloneqq \int_{\R^n} \scal{f(x)}{g(x)}_{\mathcal H} \, \diff x, \quad (f,g) \in \L\prth{\R^n,\mathcal H}^2.
\]
For an operator $A$ acting on $\L\prth{\R^n,\mathcal H}$, we denote by $A^\dagger$ its formal adjoint with respect to this inner product.
For $\ep > 0$ and $f \in \mathrm L^1\prth{\R^2, \C^2}$, we define the semiclassical Fourier transform by
\[
    \mathcal F_\ep f(\xi) \coloneqq \frac{1}{2\pi\ep} \int_{\R^2} e^{-\i \frac{x\cdot \xi}{\ep}} f(x) \diff x, \quad \xi \in \R^2,
\]
and we denote by $\widehat f$ the usual Fourier transform of $f$, related to $\mathcal F_1 f$ by $\widehat f \coloneqq 2\pi \mathcal F_1 f$.

\section{Two-scale analysis approach} \label{sec:Two_scale_analysis}

The coordinates $(s,y)$ are used to obtain a normal form of $\H$, as detailed in Section~\ref{sec:normal geodesic coordinates}.
In Section~\ref{sec:Two-scale Wigner observables}, we separate and highlight two-scale Wigner observables corresponding to the different concentration regimes.
In Section~\ref{sec:Two-scale semiclassical measures}, we prove the existence of two-scale Wigner measures and describe them for arbitrary families.
Finally, in Section~\ref{sec:rescaling}, a rescaling is applied to take into account the two-scale Wigner approach for the Hamiltonian.

\subsection{Straightening the edge} \label{sec:normal geodesic coordinates}

In this section, we consider $\prth{\psi^\ep_t}_{\ep > 0}$ a family of solutions to~\eqref{eq:Dirac} with normalized initial condition $\prth{\psi^\ep_0}_{\ep > 0}$.

\paragraph{Normal form operator.} We assume Assumption~\ref{ass:transversality} and Assumption~\ref{ass:tubularneighborhood}.
We work in the normal geodesic coordinates denoted $(s,y)$ introduced in~\eqref{eq:phi}.
Using the map $\Phi$, we transform $\H$ to a normal form, following the approach of~\cite{Ourmieres}.

\begin{proposition}[Hamiltonian Normal form near the interface] \label{pro:normal form} There exists a unitary operator $\mathrm U : \L(\Omega,\C^2) \to \L(\R\times \mathrm I,\C^2)$ such that $\mathrm U \H \mathrm U^{-1} = \H^\E$ where the operator $\H^\E$ acts on $\L(\R \times \mathrm I, \C^2)$ and is defined by
\[
    \H^\E \coloneqq \m\Big(\Phi(s,y)\Big) \sigma_1 + \ep \D_y \sigma_2 - \frac{\ep}{1 + y \kappa(s)} \D_s \sigma_3 + \frac{\i\ep y \kappa'(s)}{2(1+y\kappa(s))^2} \sigma_3.
\]
\end{proposition}

\begin{rem*} The following properties hold.
\begin{enumerate}
	\item For all $(s,y) \in \R\times \mathrm I$, we have $\abs{1+y\kappa(s)} > 1/2$; therefore, all the quantities introduced above are well-defined.
	\item The operator $\mathrm U$ is a matrix-valued Fourier Integral Operator.
	It incorporates the pull-back associated with the change of variables $\Phi$ as well as multiplication by an suitable $s$-dependent matrix, which combines a rotation and a multiplication by a $s$-dependent phase factor.
	\item In what follows, slight abuse of notation, we shall apply the operator $\mathrm U$ to functions in $\L(\R^2, \C^2)$ that are compactly supported in $\Omega$.
\end{enumerate}
\end{rem*}

\begin{proof}[Proof of Proposition~\ref{pro:normal form}]
Let us define $\mathrm U$ as follows
\[
    \mathrm U u (s,y) \coloneqq \sqrt{\frac{1 + y \kappa(s)}{2}} \begin{pmatrix} e^{ \i\frac{\theta(s)}{2}} & -e^{-\i\frac{\theta(s)}{2}} \\ e^{ \i\frac{\theta(s)}{2}} & e^{-\i\frac{\theta(s)}{2}} \end{pmatrix} u\Big(\Phi(s,y)\Big), \quad u \in \L(\Omega, \C^2), \ (s,y) \in \R \times \mathrm I.
\]
To facilitate computations, we split $\mathrm U$ into the composition of two operators.
Define
\[
    \mathrm V_\Phi u(s,y) \coloneqq \sqrt{1 + y \kappa(s)} u\Big(\Phi(s,y)\Big), \quad \quad u \in \L(\Omega, \C^2).
\]
The operator $\mathrm V_\Phi$ is unitary from $\L\prth{\Omega,\C^2}$ to $\L\prth{\R\times \mathrm I,\C^2}$, and the operator $\mathrm U_\theta$, defined in~\eqref{eq:U_theta}, is unitary on $\L\prth{\R\times \mathrm I,\C^2}$.
Notably, $\mathrm U = \mathrm U_{\theta(s)} \mathrm V_\Phi$.
Moreover, we have
\begin{align*}
    \mathrm V_\Phi \H \mathrm V_\Phi^{-1} = & \ \m\Big(\Phi(s,y)\Big) \sigma_3 - \frac{\ep}{1+y\kappa(s)} (\mathbf{t}'\odot\sigma) \D_s + \ep (\mathbf{n}\odot\sigma) \D_y \\
    & \ + \ \frac{\i\ep y\kappa'(s)}{2(1+y\kappa(s))^2} (\mathbf{t}' \odot\sigma) + \frac{\i\ep\kappa(s)}{2(1+y\kappa(s))} (\mathbf{n} \odot\sigma),
\end{align*}
where $x \odot \sigma \coloneqq x_1 \sigma_1 + x_2 \sigma_2$ for $x \in \C^2$ and $\sigma_1, \sigma_2$ are the Pauli matrices.
Finally, the following identities hold.
\begin{equation*}
\begin{split}
    \mathrm U_{\theta(s)} \sigma_3 \mathrm U_{\theta(s)}^{-1} = \sigma_1, \quad & \mathrm U_{\theta(s)} (\mathbf{t}' \odot \sigma) \mathrm U_{\theta(s)}^{-1} = \sigma_3, \\
    \mathrm U_{\theta(s)} (\mathbf{n} \odot \sigma) \mathrm U_{\theta(s)}^{-1} = \sigma_2, \quad & \mathrm U_{\theta(s)} (\mathbf{t}'\odot\sigma) \left( \D_s \mathrm U_{\theta(s)}^{-1} \right) = -\frac{\kappa(s)}{2} \i \sigma_2.
\end{split}
\end{equation*}
\end{proof}

Let us define $\displaystyle{\mathrm I_{\delta_0} \coloneqq \left( -\frac{\delta_0}{2|| \kappa ||_{\infty}} , \frac{\delta_0}{2|| \kappa ||_{\infty}}\right)}$ for $\delta_0 \in \left(0,1\right]$ and let us consider the map 
\[
    \begin{matrix} \Phi_{\delta_0} : & \R \times \mathrm I_{\delta_0} & \longrightarrow & \R^2 \\ & (s,y) & \longmapsto & \gamma(s) + y \nu(s) \end{matrix}
\]
A consequence of Assumption~\ref{ass:tubularneighborhood} is that for all $\delta_0 \in \left(0,1\right]$, the function $\Phi_{\delta_0}$ is a diffeomorphism from $\R \times \mathrm I_{\delta_0}$ to
\begin{align*}
    \Omega_{\delta_0} & \coloneqq \left\{ \gamma(s) + y \nu(s) | \ (s,y) \in \R \times \mathrm I_{\delta_0} \right\}, \\
    & = \left\{ x \in \R^2 \Big| \ \diff\prth{x,\E} < \frac{\delta_0}{2|| \kappa ||_{\infty}} \right\}.
\end{align*}
Let $\chi \in \mathscr C_c^\infty\prth{\R^2,[0,1]}$ be a cutoff function such that $\chi \equiv 1$ for $\diff(x,\E)< \frac{1}{4|| \kappa ||_{\infty}}$ and $\chi \equiv 0$ for $\diff(x,\E) > \frac{1}{2|| \kappa ||_{\infty}}$.
We then define
\begin{equation} \label{def:varphi}
    \varphi^\ep_t \coloneqq \mathrm U \chi \psi^\ep_t.
\end{equation}

\begin{lemma}[Dirac equation in normal geodesic coordinates] \label{lem:normal form} The family $\prth{\varphi^\ep_t}_{\ep > 0}$ satisfies, on $\R_t \times \R_s \times \mathrm I_{1/2}$, the following equation
\begin{equation} \label{eq:Dirac_change_variable}
    \ep \D_t \varphi^\ep_t + \mathrm H^\E_\ep \varphi^\ep_t = 0.
\end{equation}
\end{lemma}

\begin{proof} We apply $\mathrm U \chi$ to equation~\eqref{eq:Dirac}, obtaining
\[
    \ep \D_t \varphi_t^\ep + \mathrm H ^ \mathrm E _\ep \varphi_t^\ep = \mathrm U \crch{\H, \chi} \psi_t^\ep.
\]
A direct computation gives
\[
    \crch{\H, \chi} = -\i\ep \matp{0 & \partial_1 \chi - \i\partial_2 \chi \\ \partial_1 \chi + \i\partial_2 \chi & 0},
\]
where for $j \in \ens{1,2}, \partial_j \chi \in \mathscr C _c ^\infty (\R^2,\C)$ is supported in $\ens{x \in \R^2 \ \Big | \  \frac{1}{4|| \kappa ||_{\infty}} \leqslant \diff(x,\E) \leqslant \frac{1}{2|| \kappa ||_{\infty}} }$.
It follows that,
\[
    \Big( \mathrm U \crch{\H, \chi} \psi_t^\ep \Big) (s,y) = 0, \quad (t,s,y) \in \R_t \times \R_s \times \mathrm I_{1/2}.
\]
\end{proof}

\paragraph{Pull-back of semiclassical measure.}
We establish the connection between the semiclassical measure of the family $\prth{\psi_t^\ep}_{\ep > 0}$ and the two-scale semiclassical measure of the family $\prth{\varphi_t^\ep}_{\ep > 0}$ defined in~\eqref{def:varphi} via the change of variables introduced in Proposition~\ref{pro:normal form}.

\begin{lemma}[Semiclassical measure and normal geodesic coordinates] \label{lem: normal form measure} All semiclassical measures $\mu_t$ of the family $(\psi^\ep_t)_{\ep > 0}$ associated with the sequence $(\ep_k)_{k \in \N}$ satisfy
\[
    \diff\mu_t \diff t \mathds{1}_\mathcal C = \Phi^* \Big( \mathrm U_{\theta(s)}^\dagger \diff\rho_t \mathrm U_{\theta(s)} \Big) \diff t\mathds{1}_\mathcal C
\]
where $\rho_t$ is the semiclassical measure of the family $\prth{\varphi_t^\ep}_{\ep > 0}$ on $\R_s \times \mathrm I_{1/2} \times \R_{\sigma, \eta}^2 \times \R_t$ associated with the sequence $\prth{\ep_k}_{k \in \N}$.
\end{lemma}

\begin{proof} Let $a \in \mathscr C _c ^\infty \prth{\R^4, \C^{2,2}}$ with $\supp(a) \subset \Omega_{1/2} \times \R^2$.
Then, for all $N \in \mathbb{N}$, we have
\begin{align*}
    \scal{\Op \prth{a} \psi^\ep_t}{\psi^\ep_t}_{\L\prth{\R^2,\C^2}} & = \scal{\Op \prth{a} \chi \psi^\ep_t}{\chi\psi^\ep_t}_{\L\prth{\Omega_{1/2},\C^2}} + \gO{\ep^N}{} \\
    & = \scal{\mathrm U \Op \prth{a} \mathrm U ^\dagger \varphi^\ep_t}{\varphi^\ep_t}_{\L\prth{\R \times \mathrm I_{1/2},\C^2}} + \gO{\ep^N}{} .
\end{align*}
Since $\rm U = U_{\theta(s)} V_\Phi$, with $\rm U_{\theta(s)}$ depending only on $s$ and $\rm V_\Phi$ a change of variables, standard results in semiclassical pseudodifferential calculus (see~\cite[Theorem 9.3]{Zworsky}) give 
\begin{align*}
    \scal{\Op \prth{a} \psi^\ep_t}{\psi^\ep_t}
    = & \ \scal{\mathrm U_{\theta(s)} \Op \Big(\mathrm V_\Phi a \mathrm V_\Phi^\dagger \Big) \mathrm U_{\theta(s)}^\dagger \varphi^\ep_t}{\varphi^\ep_t}_{\L\prth{\R \times \mathrm I_{1/2},\C^2}} + \gO{\ep}{}
    \\
    = & \ \scal{\Op \prth{ \mathrm U_{\theta(s)} a \mathrm U_{\theta(s)}^\dagger \left( \Phi (s,y), {}^t \nabla \left( \Phi^{-1} \right) \Big( \Phi (s,y) \Big) \crch{\sigma,\eta} \right)} \varphi^\ep_t}{\varphi^\ep_t}_{\L\prth{\R \times \mathrm I_{1/2},\C^2}} \\
    & + \gO{\ep}{}.
\end{align*}
Assuming the sequence $\prth{\ep}_{\ep > 0}$ realizes the semiclassical measure $\mu_t$ of $(\psi^\ep_t)_{\ep > 0}$, passing to the limit $\ep \to 0$ gives
\begin{align*}
    \scal{\Op \prth{a} \psi^\ep_t}{\psi^\ep_t} 
    & \underset{\ep \to 0}{\longrightarrow} \int_{\R \times \mathrm I_{1/2} \times \R^2} \tr \prth{  \mathrm U_{\theta(s)} a \mathrm U_{\theta(s)}^\dagger \left( \Phi (s,y), {}^t \nabla \left( \Phi^{-1} \right) \Big( \Phi (s,y) \Big) \crch{\sigma, \eta} \right) \diff \rho_t (s,y,\sigma, \eta) } \\
    & = \int_{\R \times \mathrm I_{1/2} \times \R^2} \tr \prth{ a \left( \Phi (s,y), {}^t \nabla \left( \Phi^{-1} \right) \Big( \Phi (s,y) \Big) \crch{\sigma, \eta} \right) \Big( \mathrm U_{\theta(s)}^\dagger \diff \rho_t (s,y,\sigma, \eta) \mathrm U_{\theta(s)} \Big) } \\
    & = \int_{\Omega_{1/2} \times \R^2} \tr \prth{ a (x,\xi) \Phi^*\Big( \mathrm U_{\theta(s)}^\dagger \diff \rho_t \mathrm U_{\theta(s)} \Big) (x,\xi) }.
\end{align*}
Finally, restricting the measure $\rho_t$ to the curve $\mathcal C$ gives the semiclassical measure associated with $(\psi^\ep_t)_{\ep > 0}$ near $\mathcal C$.
\end{proof}

\subsection{Two-scale analysis of concentrations} \label{sec:Two-scale Wigner observables}

In this section, we will consider two-scale semiclassical observables and quantization defined in~\eqref{eq:secondquantization} at the scale $\sqrt \ep$.
We first remark that for all $a \in \mathcal{A}$,
\begin{align}
    \Opp(a) & = \Op\prth{a\prth{s,y,\sigma, \eta, \frac{(y,\sigma,\eta)}{\sqrt\ep}}} \nonumber \\
    & = \Lambda_{\sqrt\ep}^\dagger \mathrm{Op}_1 \Big( a \left( s,\sqrt\ep y , \ep \sigma, \sqrt\ep \eta , y, \sqrt\ep \sigma, \eta \right) \Big) \Lambda_{\sqrt\ep}, \ \prth{ \text{semiclassical in } \sqrt\ep } \label{eq:double_scale_semiclassical}
\end{align}
with
\begin{equation} \label{eq:operator_to_semiclassical}
    \Lambda_{\sqrt\ep} [f] (s,y) = \ep^{1/4} f \left( s, \sqrt\ep y  \right), \quad f \in \L\left(\R^2, \C^2\right), \ (s,y) \in \R^2.
\end{equation}
Notably, $\Lambda_{\sqrt\ep}$ is a unitary operator on $\L\prth{\R^2,\C^2}$.
Moreover, the quantization map $\Opp$ associates to $a$ a bounded operator on $\L\prth{\R^2,\C^2}$ according to the Calder\'on-Vaillancourt Theorem~\cite{CV} applied to expression~\eqref{eq:double_scale_semiclassical}.

\bigskip

These observables allow us to differentiate three regions of the phase space $\R_s \times \mathrm I_{1/2} \times \R^2_{\sigma,\eta}$.
Let us explain this fact.
Let $a \in \mathcal{A}$ and $\chi \in \mathscr C^\infty_c(\R,[0,1])$ supported in $[-1,1]$ with $\chi \equiv 1$ on $[-1/2,1/2]$.
Let $\delta > 0$, $R > R_0$ with $R_0$ associated with $a$ and $\ep$ sufficiently small enough such that
\begin{equation} \label{eq:relation_ep_delta}
    2R\sqrt\ep < \delta.
\end{equation}
Let us denote $\mathbf{x} \coloneqq (s,y), \widetilde{\xi} \coloneqq (\sigma, \eta)$ and $z \coloneqq \prth{y, \sigma, \eta}$ for simplicity.
We now split $a$ into three observables, 
\begin{align} \label{eq:observable_decomposition}
    a \left( \mathbf{x} , \widetilde{\xi} , z \right) = & \ a \left( \mathbf{x} , \widetilde{\xi} , z \right) \left( 1 - \chi\left( \frac{\sqrt\ep \abs z}{\delta} \right) \right) & \left(\eqqcolon a^\delta_\ep \left( \mathbf{x},\widetilde{\xi},z \right)\right) \\
    & + \ a \left( \mathbf{x} , \widetilde{\xi} , z \right) \chi\left( \frac{\sqrt\ep \abs z}{\delta} \right) \left( 1 - \chi\left( \frac{\abs z}{R} \right) \right) & \left(\eqqcolon a_{\ep,\delta}^R \left( \mathbf{x},\widetilde{\xi},z \right)\right) \label{eq:observable_delta} \\
    & + \ a \left( \mathbf{x} , \widetilde{\xi} , z \right) \chi\left( \frac{\abs z}{R} \right) & \left(\eqqcolon a_{\ep,R} \left( \mathbf{x},\widetilde{\xi},z \right)\right). \label{eq:observable_ep}
\end{align}
Let us define $\chi_R \coloneqq \chi\left( \frac{\abs \cdot}{R} \right)$ and $\chi^R \coloneqq 1 - \chi_R$.
The symbol $a_\ep^\delta$ through the second-scale quantization is supported outside $\mathcal C$, at a distance $\delta$ of $\mathcal C$.
Because of~\eqref{eq:relation_ep_delta}, $\Opp \left( a_\ep^\delta \right)$ is a semiclassical operator which satisfies
\[
    \Opp \left( a_\ep^\delta \right) = \mathrm{Op}_1 \left( a_\infty \left( s,y , \ep \sigma, \ep \eta , \frac{(y, \ep \sigma, \ep \eta)}{|(y, \ep \sigma, \ep \eta)|} \right) \left( 1 - \chi \left( \frac{|(y, \ep \sigma, \ep \eta)|}{\delta} \right) \right) \right).
\]
So, we can apply the semiclassical theory recalled in Appendix~\ref{app:separable} to the observable $a_\ep^\delta$.
We are particularly interested in the two regions covered by $a_{\ep,\delta}^R$ and $a_{\ep, R}$ as they describe the measure above the curve $\mathcal C$.

\bigskip

In the second-scale quantization, the symbol $a_{\ep, R}$ is supported at a distance of order $\sqrt\ep$ of $\mathcal C$ and $a_{\ep,\delta}^R$ is supported at distance larger than $\sqrt\ep$ of $\mathcal C$ but smaller than $\delta$.

\bigskip

For the observable $a_{\ep,\delta}^R$, we have
\[
    \Opp \left( a_{\ep,\delta}^R \right) = \mathrm{Op}_1 \left( a \left( s;0; \frac{y}{\sqrt\ep}, \sqrt\ep \sigma, \sqrt\ep\eta \right) \chi \left( \frac{|(y,\ep \sigma, \ep \eta)|}{\delta} \right) \left( 1 - \chi \left( \frac{|(y, \ep \sigma, \ep \eta)|}{R\sqrt\ep} \right) \right) \right) + \gO{\delta}{}.
\]
But we choose $R > R_0$ then
\[
    \Opp \left( a_{\ep,\delta}^R \right) = \mathrm{Op}_1 \left( a_\infty \left( s;0; \frac{(y,\ep\sigma,\ep\eta)}{\abs{(y,\ep\sigma,\ep\eta)}} \right) \chi \left( \frac{|(y,\ep \sigma, \ep \eta)|}{\delta} \right) \left( 1 - \chi \left( \frac{|(y, \ep \sigma, \ep \eta)|}{R\sqrt\ep} \right) \right) \right) + \gO{\delta}{}.
\]
By using~\eqref{eq:double_scale_semiclassical} on $a_{\ep, R}$, we have
\begin{align}
    \Opp \left( a_{\ep, R} \right) & = \mathrm{Op}_1 \left( a(s, 0, 0, 0, \frac{y}{\sqrt\ep}, \sqrt\ep \sigma, \sqrt\ep\eta) \chi \left( \frac{|(y,\ep \sigma, \ep\eta)|}{R\sqrt\ep} \right)\right) \ + \ \gO{\sqrt\ep}{}, \nonumber \\
    & = \Lambda_{\sqrt\ep}^\dagger \mathrm{Op}_1 \left( a(s,0, 0, 0, y, \sqrt\ep \sigma, \eta) \chi \left( \frac{|(y,\sqrt\ep \sigma, \eta)|}{R} \right)\right) \Lambda_{\sqrt\ep} \ + \ \gO{\sqrt\ep}{}. \label{eq:quantization_ep}
\end{align}
So we will consider the separable Hilbert space $\mathcal H \coloneqq \L\left(\R_y, \C^2\right)$ and the setting of Appendix~\ref{app:separable}: the functional space $\mathscr C_c^\infty \prth{\R^2, \mathcal L \prth{\mathcal H}}$, and the two-scale quantization defined for $b \in \mathscr C_c^\infty \prth{\R^2, \mathcal L \prth{\mathcal H}}$,
\begin{equation} \label{eq:second_quantization}
    \Opd (b) \coloneqq \mathrm{Op}_1 \Big( b\prth{s,\sqrt\ep\sigma} \Big).
\end{equation}
Let us denote, for all $(s,\sigma) \in \R^2$, the bounded operator on $\mathcal H$ obtained by the Weyl quantization of $(y,\eta)\mapsto a(s,0,0,0,y,\sigma,\eta)$ by $a^W (s,0,0,0, y,\sigma,\D_y)$. Then, for the operator $b(s,\sigma)$ defined by for all $(s,\sigma) \in \R^2$, $b(s,\sigma) \coloneqq a^W(s,0, 0, 0, y, \sigma, \D_y) \chi \left( \frac{|(y, \sigma, \D_y)|}{R} \right)$, we have
\begin{equation} \label{eq:second_quantization_semiclassical}
    \Opp \left( a_{\ep, R} \right) = \Lambda_{\sqrt\ep}^\dagger \Opd (b) \Lambda_{\sqrt\ep} \ + \ \gO{\sqrt\ep}{}.
\end{equation}
In the following sections, we will take advantage of these properties of $a_{\ep,R}$ and $a_{\ep,\delta}^R$.

\subsection{Two-scale semiclassical measures} \label{sec:Two-scale semiclassical measures}

\subsubsection{Existence} \label{subsec:existence}

Let us first consider the stationary case.
We use observables from $\mathcal{A}$ to analyze the oscillations at the scales $\sqrt\ep$ and $\delta$ of uniformly bounded families $\prth{f^\ep}_{\ep > 0}$ of $\L\left(\R^2,\C^2\right)$ via $a_{\ep,R}$ and $a_{\ep,\delta}^R$, respectively.
Let us define
\[
    \mathcal I^\ep_f(a) \coloneqq \scal{\Opp (a)f^{\ep}}{f^{\ep}}_{\L(\R^2,\C^2)}, \quad a \in \mathcal{A}.
\]

\begin{theorem}[Existence of the two-scale semiclassical measure] \label{thm:existence_two_scale_measure} Let $\displaystyle{(f^\ep)_{\ep > 0}}$ uniformly bounded in $\L \prth{\R^2,\C^2}$, then, there exist a vanishing sequence of positive numbers $\displaystyle{(\ep_k)_{k \in \N}}$ and
\begin{itemize}
    \item[$\bullet$] Two measurable map $\nu_\infty$ and $\nu$, valued in the set of non-negative scalar measures on $\mathcal C \times \S^2$ and $\R^2$, respectively,
    \item[$\bullet$] A $\diff\nu_\infty$-measurable map 
    $(s,\omega)\mapsto M_\infty(s,\omega) \in \C^{2,2}$, valued in the set of Hermitian trace $1$ positive matrices,
    \item[$\bullet$] A $\diff\nu(s,\sigma)$ measurable map $(s,\sigma)\mapsto M(s,\sigma)$, valued in trace-one positive operators on $\L\left(\R, \C^2\right)$,
\end{itemize}
such that for all $a\in\mathcal{A}$ we have
\begin{align*}
    \mathcal I^{\ep_k}_f(a) \underset{k \to + \infty}{\longrightarrow} & \int_{\R^4 \backslash \mathcal{C}} \mathrm{tr} \left(  a_\infty \left( s,y,\sigma,\eta, \frac{(y,\sigma,\eta)}{|(y,\sigma,\eta)|} \right) \diff\mu(s,y,\sigma,\eta) \right) \\
    & + \ \int_{\R \times \S^2} \mathrm{tr} \Big( a_\infty (s,0,0,0, \omega) M_\infty(s,\omega) \Big) \diff\nu_\infty(s,\omega) \\
    & + \ \int_{\R^2} \tr \Big( a^W (s,0,0,0,y,\sigma,\D_y) M(s,\sigma) \Big) \diff\nu(s,\sigma),
\end{align*}
where $\mu$ is the semiclassical measure of $\displaystyle{\prth{f^\ep}_{\ep > 0}}$ associated with the sequence $\displaystyle{(\ep_k)_{k \in \N}}$.
\end{theorem}

\begin{proof} Let $\displaystyle{\prth{f^\ep}_{\ep > 0}}$ uniformly bounded in $\L \prth{\R^2,\C^2}$.
For simplicity, we assume that the sequence $\prth{\ep}_{\ep > 0}$ realizes the semiclassical measure $\mu$ of $\displaystyle{\prth{f^\ep}_{\ep > 0}}$.
We use the decomposition~\eqref{eq:observable_decomposition}.
First, we analyze the concentration of $a^\delta$.
\begin{itemize}
\item[$\bullet$] Construction of $\mu$ outside $\mathcal C$
\end{itemize}
\noi For $a \in \mathcal{A}$ such that $a = 0$ in a neighborhood of $\mathcal C$ and for $\ep$ small enough, we have
\[
	\Opp \left( a \right) = \Op \left( a \left( s,y ; \sigma, \eta ; \frac{y}{\sqrt\ep}, \frac{\sigma}{\sqrt\ep}, \frac{\eta}{\sqrt\ep} \right) \right) = \mathrm{Op}_1 \left( a_\infty \left( s,y , \ep \sigma, \ep \eta , \frac{(y, \ep \sigma, \ep \eta)}{|(y, \ep \sigma, \ep \eta)|} \right) \right).
\]
After extraction of a subsequence that we still denote by $(\ep)_{\ep > 0}$, we have
\[
	\limsup_{\delta \to 0} \limsup_{R \to +\infty} \limsup_{\ep \to 0} \mathcal I^\ep_f \prth{a^\delta} = \int_{\R^4 \backslash \mathcal{C}} \mathrm{tr} \left(  a_\infty \left( s,y,\sigma,\eta, \frac{(y,\sigma,\eta)}{|(y,\sigma,\eta)|} \right) \diff\mu(s,y,\sigma,\eta) \right).
\]
Then we construct the pair $M\diff\nu$.
\begin{itemize}
\item[$\bullet$] Construction of $M\diff\nu$
\end{itemize}
\noi Let $a \in \mathscr C_c^\infty \prth{\R^7,\C^{2,2}}$.
Then, the function $b$, given by the relation~\eqref{eq:second_quantization_semiclassical}, is valued in the set of compact operators on $\mathcal H$.
We can see $\displaystyle{\prth{f^\ep}_{\ep > 0}}$ as a bounded family of $\L\prth{ \R_s, \mathcal H }$.
We are thus doing classical calculus on $\L\prth{ \R_s, \mathcal H }$. We apply Lemma~\ref{lem:mesure opératorielle}, then there exists a vanishing sequence of positive numbers $\prth{\ep_k}_{k \in \N}$ such that for all $a \in \mathscr C_c^\infty \prth{\R^7,\C^{2,2}}$,
\[
    \lim_{k \to +\infty} \scal{\Opp \left( a \right) f^{\ep_k}}{f^{\ep_k}}_{\L(\R^2,\C^2)} = \int_{\R^2} \tr \Big( a^W (s,0;0,0;y,\sigma,\D_y) M(s,\sigma) \Big) \diff\nu(s,\sigma).
\]
By applying this property to $a_{\ep,R}$, we obtain $M\diff\nu$.
Finally we construct the pair $M_\infty\diff\nu_\infty$.
\begin{itemize}
\item[$\bullet$] Construction of $M_\infty\diff\nu_\infty$
\end{itemize}
\noi Let $a \in \mathcal{A}$.
According to the relation~\eqref{eq:quantization_ep}, we are interested in the family $\widetilde f^\ep \coloneqq \Lambda_{\sqrt\ep} f^\ep$ which is bounded in $\L\left(\R^2, \C^2\right)$.
Then, the quantity
\[
    \mathscr I^{\ep_k}_{\widetilde f} \prth{\widetilde a^R_\delta} \coloneqq \scal{\mathrm{Op}_1 \left( \widetilde a^R_\delta (s,0, 0, 0, y, \sqrt{\ep_k} \sigma, \eta) \right) \widetilde f^{\ep_k}}{\widetilde f^{\ep_k}}_{\L(\R^2,\C^2)}.
 \]
is uniformly bounded in $\ep > 0$ and $R > 1$.
But the function $a^R_\delta$ is supported at distance at least $R\sqrt\ep$ of $\mathcal C$, so
\[
    \mathscr I^{\ep_k}_{\widetilde f} \prth{\widetilde a^R_\delta} = \scal{\mathrm{Op}_1 \left( a_\infty \left( s,0, 0, 0, \frac{(y, \sqrt{\ep_k} \sigma, \eta)}{\abs{\prth{y, \sqrt{\ep_k} \sigma, \eta}}} \right) \right) \widetilde f^{\ep_k}}{\widetilde f^{\ep_k}}_{\L(\R^2,\C^2)},
\]
as soon as $R$ is large enough, depending on $\ep$.
We then deduce by a diagonal extraction argument that there exists a subsequence $\prth{\ep_n}_{n \in \N}$ of $\prth{\ep_k}_{k \in \N}$, a sequence $\prth{R_n}_{n \in \N}$ and a sequence $\prth{\delta_n}_{n \in \N}$ such that for all $a \in \mathcal{A}$,
\[
    \mathscr I^{\ep_n}_{\widetilde f} \prth{\widetilde a^{R_n}_{\delta_n}} \underset{n \to +\infty}{\longrightarrow} \mathscr I \prth{a_\infty}.
\]
It remains to prove that $a_\infty \mapsto \mathscr I \prth{a_\infty}$ is a measure, which will allow us to define $M_\infty \diff\nu_\infty$. First, we prove that $a_\infty \mapsto \mathscr I \prth{a_\infty}$ is a non-negative matrix-valued distribution. We observe that there exists $C > 0$ such that for all $a \in \mathcal{A}$, for all $n \in \N$,
\[
    \mathscr I^{\ep_n}_{\widetilde{f}} \prth{\widetilde{a}^{R_n}_{\delta_n}} \leqslant C \mathrm{N}_d \prth{\widetilde{a}^{R_n}_{\delta_n}},
\]
with $\mathrm{N}_d$ a Schwarz semi-norm.
We have $\mathrm{N}_d \prth{\widetilde a^{R_n}_{\delta_n}} \to \mathrm{N}_d \prth{a_\infty}$ then $\mathscr I \prth{a_\infty} \leqslant C \mathrm{N}_d \prth{a_\infty}$.
Therefore $\mathscr I$ is a matrix-valued distribution.
Secondly, the operators $a \mapsto \mathrm{Op}_1 \left( a^R_\delta (s,0, 0, 0, y, \sqrt\ep \sigma, \eta) \right)$ satisfy a semiclassical calculus in the parameters $\ep$, $1/R$ and $\delta$.
We have for all $a,a_1,a_2 \in \mathcal{A}$,
\begin{itemize}
	\item[$\bullet$] Adjoint : $\mathrm{Op}_1 \left( \widetilde a^R_\delta (s,0, 0, 0, y, \sqrt\ep \sigma, \eta) \right)^\dagger = \mathrm{Op}_1 \left( \overline{\widetilde a^R_\delta} (s,0, 0, 0, y, \sqrt\ep \sigma, \eta) \right)$
	\item[$\bullet$] Symbolic Calculus : as operator of $\mathcal L \prth{ \L\prth{\R^2,\C^2} }$,
\[
    \mathrm{Op}_1 \Big( \widetilde a^R_{1,\delta} \Big) \mathrm{Op}_1 \Big( \widetilde a^R_{2,\delta} \Big) = \mathrm{Op}_1 \Big( \prth{\widetilde a_1 \widetilde a_2}^R_{1,\delta} \Big) + \gO{\sqrt\ep}{} + \gO{\frac{1}{R}}{},
\]
	\item[$\bullet$] Weak G\aa rding inequality : if $a \geqslant 0$, then for all $\widetilde \delta$, there exists $C_{\widetilde \delta} > 0$, such that for all $f \in \L\prth{\R^2,\C^2}$,
\[
    \scal{\mathrm{Op}_1 \Big( \widetilde a^R_\delta \Big) \widetilde f}{\widetilde f}_{\L\prth{\R^2,\C^2}} \geqslant -\prth{\widetilde \delta + C_{\widetilde \delta} \prth{\ep + \frac{1}{R}} }\norme{f}^2_{\L\prth{\R^2,\C^2}},
\]
\end{itemize}
Therefore $a_\infty \mapsto \mathscr I \prth{a_\infty}$ is a non-negative matrix-valued measure defined on $\R \times \S^2$ that we denote $M_\infty\diff\nu_\infty$.
\end{proof}

\subsubsection{First properties}

In this section, we first explore properties of two-scale semiclassical measures defined by Theorem~\ref{thm:existence_two_scale_measure}.
A direct consequence of the proof is the separated construction of limits as follows.

\begin{remark}[Link with the limit of the previous observable] \label{rem:link_measure_observable} Let $\displaystyle{(f^\ep)_{\ep > 0}}$ uniformly bounded in $\L \prth{\R^2,\C^2}$ and a vanishing sequence of positive numbers $\displaystyle{(\ep_k)_{k \in \N}}$ given by Theorem~\ref{thm:existence_two_scale_measure}, then, with the same notations, for all $a \in \mathcal{A}$,
\begin{align*}
	\limsup_{\delta \to 0} \limsup_{R \to +\infty} \lim_{k \to +\infty} \mathcal I^{\ep_k}_f\left( a_{\ep_k}^\delta \right) = & \int_{\R^4 \backslash \mathcal{C}} \mathrm{tr} \left( a_\infty \left( s,y,\sigma,\eta, \frac{(y,\sigma,\eta)}{|(y,\sigma,\eta)|} \right) \diff\mu(s,y,\sigma,\eta) \right), \\
	\limsup_{\delta \to 0} \limsup_{R \to +\infty} \lim_{k \to +\infty} \mathcal I^{\ep_k}_f \left( a_{\ep_k, \delta}^R \right) = & \int_{\R \times \S^2} \mathrm{tr} \Big( a_\infty (s,0,0,0, \omega) M_\infty(s,\omega) \Big) \diff\nu_\infty(s,\omega), \\
	\limsup_{\delta \to 0} \limsup_{R \to +\infty} \lim_{k \to +\infty} \mathcal I^{\ep_k}_f \left( a_{\ep_k,R} \right) = & \int_{\R^2} \tr \Big( a^W (s,0,0,0,y,\sigma,\D_y) M(s,\sigma) \Big) \diff\nu(s,\sigma).
\end{align*}
\end{remark}

\begin{exe*} Let us consider,
\[
    \phi_\ep(s,y) = \frac{1}{\ep^{1/4}} u_\ep\prth{s} f\prth{\frac{y}{\sqrt\ep}} \vec V, \quad (s,y) \in \R^2,
\]
with $\vec V \in \C^2$, $(f,u_\ep) \in \mathscr S (\R,\C)^2$ and the knowledge of $\rho$ the semiclassical measure of $(u_\ep)_{\ep > 0}$ at the scale $\sqrt\ep$, then
\[
    \scal{\Opp(a)\phi_\ep}{\phi_\ep} \underset{\ep \to 0}{\longrightarrow} \int_{\R^2} \tr \Big( a^W (s,0,0,0,y,\sigma,\D_y) M(s,\sigma) \Big) \diff\nu(s,\sigma), \quad a \in \mathscr{C}^\infty_c (\R^7,\C^{2,2}),
\]
where
\[
    M(s,\sigma)\diff\nu (s,\sigma) = \frac{1}{2\pi} \Pi_{f\vec V} \diff\rho(s,\sigma),
\]
with $\Pi_{f\vec V}$ the projector on $\L\prth{\R,\C^2}$ over the function $f\vec V$.
Notably, if $u_\ep$ is a wave packet, which means
\[
    u_\ep\prth{s} \coloneqq \frac{1}{\ep^{1/4}} u\prth{\frac{s}{\sqrt\ep}}, \quad u \in \mathscr S (\R,\C), \quad s \in \R,
\]
then the semiclassical measure of $(u_\ep)_{\ep > 0}$ at the scale $\sqrt\ep$ is $\rho(\sigma) = \left| \mathcal F u (\sigma) \right|^2$ (according to the case $\alpha = 1$ of Example~\ref{example1} but with a $\sqrt\ep$-Weyl quantization) so
\[
    M(s,\sigma)\diff\nu (s,\sigma) =  \left| \mathcal F u (\sigma) \right|^2 \Pi_{f\vec V} \frac{\diff\sigma}{2\pi}.
\]
\end{exe*}

Another consequence is the application to observables in $\mathscr C_c^\infty \prth{\R^4,\C^{2,2}}$.

\begin{remark} \label{rem:measure_observable_phase_space} Let $\displaystyle{(f^\ep)_{\ep > 0}}$ uniformly bounded in $\L \prth{\R^2,\C^2}$ and a vanishing sequence of positive numbers $\displaystyle{(\ep_k)_{k \in \N}}$ given by Theorem~\ref{thm:existence_two_scale_measure}. According to the notations of Theorem~\ref{thm:existence_two_scale_measure}, the two-scale semiclassical measures $\displaystyle{\prth{M_\infty\diff\nu_\infty, M\diff\nu}}$ satisfy for all $a \in \mathscr C_c^\infty \prth{\R^4,\C^{2,2}}$,
\begin{align*}
    \mathcal I^{\ep_k}_f(a) \underset{k \to + \infty}{\longrightarrow} & \int_{\R^4 \backslash \mathcal{C}} \mathrm{tr} \Big(  a \left( s,y,\sigma,\eta \right) \diff\mu(s,y,\sigma,\eta) \Big) + \int_{\R \times \S^2} \mathrm{tr} \Big( a (s,0,0,0) M_\infty(s,\omega) \Big) \diff\nu_\infty(s,\omega) \\
    & + \ \int_{\R^2} \tr \Big( a(s,0,0,0) M(s,\sigma) \Big) \diff\nu(s,\sigma).
\end{align*}
\end{remark}

Even though we cannot describe the evolution of the two-scale semiclassical measures at infinity $M_\infty\diff\nu_\infty$, we have a vanishing condition.

\begin{corollary}[A condition to have $M_\infty\diff\nu_\infty=0$] \label{cor:condition_M_infty_vanishing} Let $\displaystyle{\prth{f^\ep}_{\ep > 0}}$ uniformly bounded in $\L \prth{\R^2,\C^2}$ such that there exists $\chi$ a cut off function with $\chi \in \mathscr{C}^\infty_c(\R,\R)$ supported in $[-1,1]$ and $\chi \equiv 1$ on $[-1/2,1/2]$ such that
\[
    \limsup_{\ep \to 0} \norme{\mathrm{Op}_{1} \prth{ 1 - \chi \left( \frac{|(y, \ep \sigma, \ep \eta)|}{R\sqrt\ep} \right) } f^\ep}_{\L\prth{\R^2,\C^2}} \underset{R \to +\infty}{\longrightarrow} 0.
\]
Then, for all two-scale semiclassical measure $M_\infty\diff\nu_\infty$ at infinity of $\displaystyle{\prth{f^\ep}_{\ep > 0}}$,
\[
    M_\infty\diff\nu_\infty = 0.
\]
\end{corollary}

\begin{proof}
Let $\displaystyle{\prth{f^\ep}_{\ep > 0}}$ uniformly bounded in $\L \prth{\R^2,\C^2}$ and $\chi$ a cut off function, $\chi \in \mathscr{C}^\infty_c(\R,[0,1])$ supported in $[-1,1]$ and $\chi \equiv 1$ on $[-1/2,1/2]$ such that
\[
    \limsup_{\ep \to 0} \norme{\mathrm{Op}_{1} \prth{ 1 - \chi \left( \frac{|(y, \ep \sigma, \ep \eta)|}{R\sqrt\ep} \right) } f^\ep}_{\L\prth{\R^2,\C^2}} \underset{R \to +\infty}{\longrightarrow} 0.
\]
For all $a\in\mathcal{A}$, we consider, 
\begin{align*}
    \alpha(s,y,\sigma,\eta) & \coloneqq a \left( s, y, \ep\sigma, \ep\eta, \frac{y}{\sqrt\ep}, \sqrt\ep \sigma, \sqrt\ep\eta \right) \chi \left( \frac{|(y,\ep \sigma, \ep \eta)|}{\delta} \right), \\
    \beta(y,\sigma,\eta) & \coloneqq 1 - \chi \left( \frac{|(y, \ep \sigma, \ep \eta)|}{R\sqrt\ep} \right).
\end{align*}
According to pseudodifferential calculus and of the inequality~\eqref{eq:relation_ep_delta}, we fix $\delta$ and let $\ep$ goes to $0$ then $R$ goes to $+\infty$, so
\[
    \mathrm{Op}_1 \Big( \alpha(s,y,\sigma,\eta) \beta(y,\sigma,\eta) \Big) = \mathrm{Op}_1 \Big( \alpha(s,y,\sigma,\eta) \Big) \mathrm{Op}_1 \Big( \beta(y,\sigma,\eta) \Big) + \gO{\sqrt\ep}{} + \gO{\frac{1}{R}}{},
\]
where remainders are taken in $\mathcal{L}\prth{\L\prth{\R^2,\C^2}}$.
Moreover,
\[
    \mathrm{Op}_1 \Big( \alpha(s,y,\sigma,\eta) \Big) = \Lambda_{\sqrt\ep}^\dagger \mathrm{Op}_1 \left( a \left( s, \sqrt\ep y, \ep\sigma, \sqrt\ep\eta, y, \sqrt\ep \sigma, \eta \right) \chi \left( \frac{\sqrt\ep|(y,\sqrt\ep \sigma, \eta)|}{\delta} \right) \right) \Lambda_{\sqrt\ep}.
\]
By application of the Calder\'on-Vaillancourt Theorem~\cite{CV},
\begin{align*}
    \norme{\mathrm{Op}_1 \Big( \alpha(s,y,\sigma,\eta) \Big)}_{\mathcal{L}\prth{\L}} & = \norme{ \mathrm{Op}_1 \left( a \left( s, \sqrt\ep y, \ep\sigma, \sqrt\ep\eta, y, \sqrt\ep \sigma, \eta \right) \chi \left( \frac{\sqrt\ep|(y,\sqrt\ep \sigma, \eta)|}{\delta} \right) \right) }_{\mathcal{L}\prth{\L}} \\
    & \leqslant C(a)
\end{align*}
where $C$ does not depend of $\ep$.
Then
\[
    \mathcal I^\ep_f \left( a_{\ep, \delta}^R \right) = \Big\langle \mathrm{Op}_{1} \Big( \alpha(s,y,\sigma,\eta) \Big) f^\ep , \mathrm{Op}_{1} \Big( \beta(y,\sigma,\eta) \Big) f^\ep \Big\rangle_{\L\prth{\R^2,\C^2}} + \gO{\sqrt\ep}{} + \gO{\frac{1}{R}}{},
\]
so,
\begin{align*}
    \left| \mathcal I^\ep_f \left( a_{\ep, \delta}^R \right) \right| \leqslant & \ \norme{\mathrm{Op}_{1} \Big( \alpha(s,y,\sigma,\eta) \Big) f^\ep}_{\L\prth{\R^2,\C^2}} \norme{\mathrm{Op}_{1} \Big( \beta(y,\sigma,\eta) \Big) f^\ep}_{\L\prth{\R^2,\C^2}} + \gO{\sqrt\ep}{} + \gO{\frac{1}{R}}{}, \\
    \leqslant & \ C(a) \norme{f^\ep}_{\L\prth{\R^2,\C^2}} \norme{\mathrm{Op}_{1} \prth{ 1 - \chi \left( \frac{|(y, \ep \sigma, \ep \eta)|}{R\sqrt\ep} \right) } f^\ep}_{\L\prth{\R^2,\C^2}} + \gO{\sqrt\ep}{} + \gO{\frac{1}{R}}{}.
\end{align*}
Then,
\[
    \limsup_{\delta \to 0} \limsup_{R \to +\infty} \limsup_{\ep \to 0} \mathcal I^\ep_f \left( a_{\ep, \delta}^R \right) = 0.
\]
We conclude by using Remark~\ref{rem:link_measure_observable}.
\end{proof}

\paragraph{Time dependency.} As in Section~\ref{subsec:wigner}, one can extend the definition of two-scale semiclassical measures to time-dependent families.
More precisely, we consider families $\displaystyle{\prth{f^\ep_t}_{\ep > 0}}$ that are uniformly bounded in $\mathrm L^\infty \prth{\R, \mathrm L^2 \prth{\R^2,\C^2}}$.
We consider
\[
    \mathcal I^\ep_f(\Xi,a) \coloneqq \int_\R \Xi(t) \scal{\Opp (a)f_t^{\ep}}{f_t^{\ep}}_{\L(\R^2,\C^2)} \diff t, \quad a\in\mathcal{A}, \quad \Xi \in \mathscr C_c^\infty(\R,\C),
\]
As before, we are interested in the limit of $\mathcal I^\ep_f(\Xi,a)$ as $\ep$ goes to $0$.
Adapting the context of ~\cite{Fermanian} to time-dependent families, one defines two-scale semiclassical measures over the curve $\mathcal C$.
Indeed, by the same argument of Section~\ref{subsec:wigner}, we can extend Theorem~\ref{thm:existence_two_scale_measure} to time-dependent families.
Therefore, as before, the limit measures will be absolutely continuous in time variable with respect to Lebesgue measure.

\subsection{Rescaling} \label{sec:rescaling}

According to the expression of the rescaling operator $\Lambda_{\sqrt\ep}$ defined in~\eqref{eq:operator_to_semiclassical}, we consider the family $\prth{u^\ep_t}_{\ep > 0}$ defined such that
\[
    \Lambda_{\sqrt\ep}^\dagger \crch{u^\ep_t} (s,y) \coloneqq \varphi_t^\ep (s,y), \quad \quad \forall (t,s,y) \in \R \times \R \times \mathrm I_{1/2}. 
\]
Although the family $\prth{u^\ep_t}_{\ep > 0}$ is defined on $\R \times \mathrm I_{\frac{1}{2\sqrt\ep}}$, we will study it on $\R \times \mathrm I_{\sqrt\ep}^\delta$ with
\[
    \mathrm I _{\sqrt\ep}^\delta \coloneqq \left( - \frac{\delta}{\sqrt\ep}, \frac{\delta}{\sqrt\ep} \right)
\]
where we consider $\delta \leqslant \frac{1}{4 \normi{\kappa}}$.

\begin{lemma}[Double scaling operator] The family $\prth{u^\ep_t}_{\ep > 0}$ satisfies on $\R\times \mathrm I _{\sqrt\ep}^\delta$,
\begin{equation} \label{eq:Dirac_rescale}
    \left( \sqrt{\ep} \D_t + y\abs{\nabla\m(\mathbf{t}(s))} \sigma_1 + \D_y \sigma_2 - \frac{\sqrt{\ep}}{1 + \sqrt{\ep} y \kappa(s)} \D_s \sigma_3 + \frac{\i\ep y \kappa'(s)}{2(1+\sqrt{\ep}y\kappa(s))^2} \sigma_3 \right) u_t^\ep = -\sqrt\ep \mathrm{R}^\ep_2 u_t^\ep,
\end{equation}
where
\[
    \big( \mathrm{R}_2^\ep u_t^\ep \big) (s,y) \coloneqq y^2 \int_0^1 \nabla^2 \m \Big(\Phi(s,\tau\sqrt\ep y)\Big) \cdot \Big( \mathbf{n}(s),\mathbf{n}(s) \Big) \diff\tau \sigma_1 u_t^\ep (s,y), \quad (s,y) \in \R \times \mathrm I _{\sqrt\ep}^\delta.
\]
\end{lemma}

\begin{proof}
On the one hand, we rewrite the equation~\eqref{eq:Dirac_change_variable} with the definition of $\prth{u^\ep_t}_{\ep > 0}$ on $\R \times \mathrm I _{\sqrt\ep}^\delta$,
\[
    \left[ \ep \D_t + \m\Big(\Phi(s, \sqrt{\ep}y)\Big) \sigma_1 + \sqrt{\ep} \D_y \sigma_2 - \frac{\ep}{1 + \sqrt{\ep} y \kappa(s)} \D_s \sigma_3 + \frac{\i\ep^{3/2} y \kappa'(s)}{2(1+\sqrt{\ep}y\kappa(s))^2} \sigma_3 \right] u^\ep_t
    = 0.
\]
On the other hand, for all $(s,y) \in \R \times \mathrm I _{\sqrt\ep}^\delta$,
\begin{align*}
	\m \left( \Phi(s, \sqrt\ep y) \right) = & \ \underbrace{\m(\mathbf{t}(s))}_{=0} + \sqrt\ep y \underbrace{\nabla \m (\mathbf{t}(s)) \cdot \mathbf{n}(s)}_{= |\nabla \m (\mathbf{t}(s))|} \\
    & + \ \ep y^2 \int_0^1 (1-\tau) \nabla^2 \m \Big(\Phi(s,\tau\sqrt\ep y)\Big) \cdot \Big( \mathbf{n}(s),\mathbf{n}(s) \Big) \diff\tau.
\end{align*}
\end{proof}

The operator $\sqrt\ep \mathrm{R}^\ep_2$ will play a role of remainder in the determination of the two-scale semiclassical measure at finite distance $M\diff\nu$.

\bigskip

For now, we make a few remarks about the two-scale operator.
The quantity $\sqrt\ep$ is in front of $\D_s$ while no power of $\ep$ multiplies $\D_y$.
It fits with the two-scale operators that we have introduced in Section~\ref{sec:main_result}.
It leads us with a semiclassical operator in the variable $s$ at scale $\sqrt\ep$ and a classical operator in the variable $y$.
Because of the second-scale of the problem, we will consider the symbol of the principal term of equation~\eqref{eq:Dirac_rescale} only with respect to the variable $s$, keeping the operator character over the variable $y$ as mentioned earlier in Section~\ref{sec:main_result}.
In fact, we had introduced the operator $\mathrm T_\E$,
the symbol of an operator acting on $\L\prth{\R_s,\L\prth{\R_y,\C^2}}$ in~\eqref{eq:second_quantization}, its eigenvalues in~\eqref{eq:lambda_scale_ep} with the associated normalized eigenvectors in~\eqref{eq:eigenvector_scale_ep} and the associated eigenprojectors in~\eqref{eq:proj_scale_ep}.

\bigskip

Let us define the symbol of $\mathrm{R}^\ep_2$ for the two-scale quantization~\eqref{eq:second_quantization} as the operator acting on $\L\left(\R_y, \C^2\right)$ such that
\[
    \mathfrak{r}^\ep_2 (s,\sigma) \coloneqq y^2 \int_0^1 \nabla^2 \m \Big(\Phi(s,\tau\sqrt\ep y)\Big) \cdot \Big( \mathbf{n}(s),\mathbf{n}(s) \Big) \diff\tau \sigma_1, \quad (s,\sigma) \in \R^2.
\]
Moreover, we define the operator acting on $\L\left(\R_y, \C^2\right)$ such that
\[
    \mathfrak{r}^\ep_3 (s,\sigma) \coloneqq \frac{ y\kappa(s) \sigma }{1 + \sqrt\ep y \kappa(s)} \sigma_3, \quad (s,\sigma) \in \R^2.
\]
We define $\mathfrak{r}^\ep \coloneqq \mathfrak{r}^\ep_2 + \mathfrak{r}^\ep_3$.
With the notation~\eqref{eq:second_quantization}, the equation~\eqref{eq:Dirac_rescale} writes
\[
    \sqrt\ep \D_t u^\ep_t + \Opd \Big( \mathrm T_\E + \sqrt\ep \mathfrak{r}^\ep \Big) u^\ep_t = 0.
\]

\section{Proof of the main result} \label{sec:proof}

In this section we describe the evolution of the two-scale semiclassical measures defined in Theorem~\ref{thm:shape_two_scale_measure} for $\prth{\varphi^\ep_t}_{\ep > 0}$ solutions to equation~\eqref{eq:Dirac_change_variable}.
In Section~\ref{sub:Gamma}, we establish the evolution of the measure at finite distance $M^t\diff\nu^t$.
In Section~\ref{sub:M}, we establish the evolution of the measure at infinity $M_\infty^t\diff\nu_\infty^t$.
In Section~\ref{sec:proof_main}, we establish Theorem~\ref{thm:main_result}.
In Section~\ref{sec:proof_applications}, we finally establish Corollaries~\ref{cor:conservation},~\ref{cor:goal_reached},~\ref{cor:linkDrouot} and~\ref{cor:below_estimate}.

\subsection{Two-scale semiclassical measures above the interface} \label{sec:section}

\subsubsection{Two-scale Wigner measure at finite distance } \label{sub:Gamma}

We prove here the first point of Theorem~\ref{thm:main_result}, which is the evolution of the measure at finite distance $M^t\diff\nu^t$.
We recall the following lemma of~\cite{GMMPErratum}.

\begin{lemma}[Poisson bracket and projector properties] \label{lem:poissonbracket} Let $\prth{\Pi_n}_{n \in \N}$ be a family of orthogonal projectors such that for all $(n,k) \in \N^2, \Pi_n\Pi_k=\delta_{k,n}$.
Then
\[
	\Pi_k \Big( \ens{\Pi_j , \Pi_\ell} - \ens{\Pi_\ell , \Pi_j} \Big) \Pi_k = 0, \quad (j,k,\ell) \in \N^3.
\]
\end{lemma}

\begin{proof}[Proof of the first point of Theorem~\ref{thm:main_result}] 
Assume Assumption~\ref{ass:transversality} and Assumption~\ref{ass:tubularneighborhood}.
For convenience, we assume that the sequence $\prth{\ep}_{\ep > 0}$ realizes the measure $M^t\diff\nu^t$ of Theorem~\ref{thm:shape_two_scale_measure}.

\medskip

Let $a \in \mathcal{A}$ and $R,\delta$ as defined in~\eqref{eq:relation_ep_delta}.
We consider the part~\eqref{eq:observable_ep} of $a$ denoted $a_{\ep, R}$ and $b(s,\sigma) \coloneqq a^W(s,0, 0, 0, y, \sigma, \D_y) \chi \left( \frac{|(y, \sigma, \D_y)|}{R} \right)$, the operator acting on $\mathcal H \coloneqq \L\left(\R_y, \C^2\right)$, obtained by Weyl quantization.
The link between the second-scale quantization of $a_{\ep, R}$~\eqref{eq:secondquantization} and the two-scale quantization of $b$~\eqref{eq:second_quantization} was established in~\eqref{eq:second_quantization_semiclassical}.
By considering $\prth{u^\ep_t}_{\ep > 0}$ solution to equation~\eqref{eq:Dirac_rescale} as a family of $\L\left(\R_s,\mathcal H\right)$, we have
\[
    \i\frac{\mathrm d}{\mathrm{dt}} \scal{ \Opd (b) u^\ep_t}{u^\ep_t} = \frac{1}{\sqrt\ep} \scal{ \Big[ \Opd (b), \Opd \Big( \mathrm T_\E \Big) \Big] u^\ep_t}{u^\ep_t} + \scal{ \Big[ \Opd (b), \Opd \big( \mathfrak{r}^\ep \big) \Big] u^\ep_t}{u^\ep_t}.
\]
Because of~\eqref{eq:poissonbracket}, we have
\begin{align*}
    \frac{\mathrm d}{\mathrm{dt}} \scal{ \Opd (b) u^\ep_t}{u^\ep_t} = & \ \frac{1}{\i\sqrt\ep} \scal{ \Opd \Big( \crch{b,\mathrm T_\E} \Big) u^\ep_t}{u^\ep_t} \\
    & + \frac{1}{\i} \scal{ \Opd \left( \frac{1}{2\i}\ens{b, \mathrm T_\E} - \frac{1}{2\i}\ens{\mathrm T_\E, b} + \crch{b,\mathfrak{r}^\ep} \right) u^\ep_t}{u^\ep_t} \\
    & + \gO{\sqrt\ep}{R}.
\end{align*}
\begin{itemize}
    \item[$\bullet$] Step 1 : Structure of $M^t$
\end{itemize}
Let $\Xi \in \mathscr C^\infty_c (\R)$, then
\begin{align*}
    \int_\R \Xi(t) \i \sqrt\ep \frac{\mathrm d}{\mathrm{dt}} \scal{\Opd (b) u^\ep_t}{u^\ep_t} \diff t = & \int_\R \Xi(t) \scal{ \Opd \Big( \crch{b,\mathrm T_\E} \Big) u^\ep_t}{u^\ep_t} \diff t \\
    & + \frac{\sqrt\ep}{2\i} \int_\R \Xi(t) \scal{ \Opd \Big( \ens{b,\mathrm T_\E} - \ens{\mathrm T_\E, b} \Big) u^\ep_t}{u^\ep_t} \diff t \\
    & + \sqrt\ep \int_\R \Xi(t) \scal{ \Opd \Big( \crch{b,\mathfrak{r}^\ep} \Big) u^\ep_t}{u^\ep_t} \diff t + \gO{\ep}{R} \\
    = & -\i\sqrt\ep \int_\R \Xi'(t) \scal{\Opd (b) u^\ep_t}{u^\ep_t} \diff t + \gO{\ep}{R}.
\end{align*}
But $b$ is compactly supported in $\ens{ \abs{(y,\sigma)} \leqslant R }$ and for all $(s,\sigma,y) \in \R \times \left( -R, R \right)^2$,
\[
    y \underbrace{\frac{ \kappa(s)}{1 + \sqrt\ep y \kappa(s)}}_{\in \mathrm L^\infty \big( \R \times \left( -R, R \right)\big)} \sigma = \gO{R^2}{}, \quad \quad y^2 \underbrace{\int_0^1 \nabla^2 \m \Big(\Phi(s,\tau\sqrt\ep y)\Big) \cdot \Big( \mathbf{n}(s),\mathbf{n}(s) \Big) \diff\tau}_{\in \mathrm L^\infty \big(\R \times \left( -R, R \right)\big)} = \gO{R^2}{}.
\]
Then
\[
    \sqrt\ep \int_\R \Xi(t) \scal{ \Opd \Big( [b,\mathfrak{r}^\ep] \Big) u^\ep_t}{u^\ep_t} \diff t = \gO{\sqrt\ep}{R}.
\]
By passing to the limit, we have, for all $\Xi \in \mathscr C^\infty_c (\R)$ and all $a\in \mathcal{A}$,
\begin{align*}
    0 & = \int_{\R^3} \Xi(t) \tr_{\L\prth{\R_y,\C^2}} \Bigg( \Big[a(s, 0 , 0 , y, \sigma, \D_y), \mathrm T_\E(s,\sigma) \Big] M^t (s,\sigma) \Bigg) \diff \nu^t (s,\sigma) \diff t \\
    & = \int_{\R^3} \Xi(t) \tr \Bigg( a(s, 0 , 0 , y, \sigma, \D_y) \Big[ \mathrm T_\E(s,\sigma), M^t (s,\sigma) \Big] \Bigg) \diff \nu^t (s,\sigma) \diff t.
\end{align*}
So, $t$ a.e. $\nu^t$ a.e. $\Big[ \mathrm T_\E(s,\sigma), M^t (s,\sigma) \Big] = 0$, which implies~\eqref{eq:measure_scale_ep}, according to the spectral decomposition of the operator $\mathrm T_\E$.

\begin{itemize}
    \item[$\bullet$] Step 2 : Equation of $\prth{\nu_n^t}_{n \in \Z}$
\end{itemize}
By Step $1$, for all $(k,\ell) \in \Z^2$, with $k \neq \ell$ and for all $\beta$ function from $\R_s\times\R_\sigma$ to $\mathcal L \prth{\L\prth{\R_y,\C^2}}$ such that $\Pi_k \beta \Pi_\ell \in \mathscr C_c^\infty \prth{\R^2, \mathcal L \prth{\mathcal H}}$, we have
\[
    \limsup_{\delta \to 0} \limsup_{R \to +\infty} \lim_{\ep \to 0} \int_\R \Xi(t) \scal{\Opd \Big( \Pi_k \beta \Pi_\ell \Big) u^\ep_t}{u^\ep_t} \diff t = 0.
\]
Then we can consider symbol of the form $\Pi_\ell \beta \Pi_\ell$.
Since the projector $\Pi_\ell$ has rank one, we have
\[
    \Pi_\ell \beta \Pi_\ell = \tr\Big( \Pi_\ell \beta \Pi_\ell \Big) \Pi_\ell,
\]
where $\tr_{\L\prth{\R_y,\C^2}}\big( \Pi_\ell \beta \Pi_\ell \big) \in \mathscr C_c^\infty \prth{\R^2, \C}$.
Therefore we will consider symbols of the form $\beta \Pi_\ell$ with $\beta \in \mathscr C_c^\infty \prth{\R^2, \C}$. Let $n \in \Z$ and $\beta \in \mathscr C_c^\infty \prth{\R^2, \C}$.
Let us define
\[
    b_n(s,\sigma) \coloneqq \beta(s,\sigma) \Pi_n(s,\sigma), \quad (s,\sigma) \in \R^2.
\]
The expression of $b_n$ implies $[b_n,\mathrm T_\E] = 0$.
Then
\begin{align*}
    \frac{\mathrm d}{\mathrm{dt}} \scal{\Opd (b_n) u^\ep_t}{u^\ep_t} = & \ \frac{1}{2} \scal{ \Opd \Big( \{\mathrm T_\E, b_n \} - \{ b_n, \mathrm T_\E\} \Big) u^\ep_t}{u^\ep_t} \\
    & + \frac{1}{\i} \scal{ \Opd \Big( [b_n,\mathfrak{r}^\ep] \Big) u^\ep_t}{u^\ep_t} + \gO{\sqrt\ep}{R}.
\end{align*}
We first show that the term with $\mathfrak{r}^\ep$ will vanish as $\ep$ goes to $0$.
By definition of $M^t\diff\nu^t$, we have
\begin{align*}
    \scal{ \Opd \big( [b_n,\mathfrak{r}^\ep] \big) u^\ep_t}{u^\ep_t} & \underset{\ep \to 0}{\longrightarrow} \int_{\R^2} \tr \Bigg( \Big[a(s, 0 , 0 , y, \sigma, D_y), \mathfrak{r}^0(s,\sigma) \Big] M^t (s,\sigma) \Bigg) \diff \nu^t (s,\sigma) \\
    & = \sum_{m \in \Z} \int_{\R^2} \tr \Bigg( \Big[ \beta(s, \sigma) \Pi_n(s,\sigma) , \mathfrak{r}^0(s,\sigma) \Big] \Pi_m (s,\sigma) \Bigg) \diff \nu^t_m (s,\sigma), \\
    & = \sum_{m \in \Z} \int_{\R^2} \beta(s, \sigma) \tr \Bigg( \Big[ \Pi_n(s,\sigma) , \mathfrak{r}^0(s,\sigma) \Big] \Pi_m (s,\sigma) \Bigg) \diff \nu^t_m (s,\sigma),
\end{align*}
where the second equality holds because of Step $1$ and the chosen shape of $a$ and the last equality holds because $\beta$ is scalar and $\mathfrak{r}^0$ is the value at zero of $\mathfrak{r}^\ep$ given by for all $(s,\sigma) \in \R^2$,
\[
    \mathfrak{r}^0(s,\sigma) = \mathfrak{r}^0_2(s,\sigma) + \mathfrak{r}^0_3(s,\sigma) = y^2 \nabla^2 \m \Big(\Phi(s,0)\Big) \cdot \Big( \mathbf{n}(s),\mathbf{n}(s) \Big) \sigma_1 + y\kappa(s) \sigma \sigma_3.
\]
Then, we recall, for all $m \in \Z$,
\begin{align*}
    \tr \Big( \big[ \Pi_n , \mathfrak{r}^0 \big] \Pi_m \Big) & = \tr \Big( \Pi_n \mathfrak{r}^0 \Pi_m \Big) - \tr \Big( \mathfrak{r}^0 \Pi_n \Pi_m \Big) = \tr \Big( \mathfrak{r}^0 \Pi_m \Pi_n \Big) - \tr \Big( \mathfrak{r}^0 \Pi_n \Pi_m \Big) \\
    & = \delta_{n=m} \left( \tr \Big( \mathfrak{r}^0 \Pi_n \Big) - \tr \Big( \mathfrak{r}^0 \Pi_n \Big) \right) = 0.
\end{align*}
Let $k \in \Z$, let us compute $ \Pi_k \{\mathrm T_\E, b_n \} \Pi_k$,
\begin{equation}
	\Pi_k \{\mathrm T_\E, b_n \} \Pi_k = \Pi_k \left\{ \sum_{m \in \Z} \lambda_m \Pi_m \, , \, \beta \Pi_n \right\} \Pi_k \label{eq:decomposition_proj_n}
\end{equation}
then,
\[
	\Pi_k \{\mathrm T_\E, b_n \} \Pi_k = \Pi_k \sum_{m \in \Z} \Bigg( \Pi_m \left\{ \lambda_m, \beta \right\} \Pi_n + \lambda_m \left\{ \Pi_m, \beta \right\} \Pi_n + \lambda_m \beta \left\{ \Pi_m, \Pi_n \right\} + \beta \Pi_m \left\{ \lambda_m, \Pi_n \right\} \Bigg) \Pi_k.
\]
Let $m \in \Z$,
\begin{align*}
	\Pi_k \Pi_m \left\{ \lambda_m, \beta \right\} \Pi_n \Pi_k & = \mathds{1}_{k=m=n} \left\{ \lambda_n, \beta \right\} \Pi_n, \\
	\Pi_k \lambda_m \left\{ \Pi_m, \beta \right\} \Pi_n \Pi_k & = \mathds{1}_{n=k} \lambda_m \Pi_n \left\{ \Pi_m, \beta \right\} \Pi_n \ = \ 0, \\
	\Pi_k \lambda_m \beta \left\{ \Pi_m, \Pi_n \right\} \Pi_k & = \lambda_m \beta \Pi_k \left\{ \Pi_m, \Pi_n \right\} \Pi_k, \\
	\Pi_k \beta \Pi_m \left\{ \lambda_m, \Pi_n \right\} \Pi_k & = \mathds{1}_{k=m}\beta \Pi_m \left\{ \lambda_m, \Pi_n \right\} \Pi_m \ = \ 0,
\end{align*}
where the two vanishing equations hold because $\lambda_n$ and $\beta$ are scalar and $\Pi_k \partial \Pi_m \Pi_k$ is always equal to zero. Moreover, according to Lemma~\ref{lem:poissonbracket},
\[
	\lambda_m \beta \Pi_k \Big( \left\{ \Pi_m, \Pi_n \right\} - \left\{ \Pi_n, \Pi_m \right\} \Big) \Pi_k = 0 \, , \, \text{for all} \ m \in \Z.
\]
Therefore,
\[
    \{\mathrm T_\E, b_n \} - \{ b_n, \mathrm T_\E\} = \Big( \left\{ \lambda_n, \beta \right\} - \left\{ \beta , \lambda_n \right\} \Big) \Pi_n.
\]
Then,
\[
    \frac{\mathrm d}{\mathrm{dt}} \scal{\Opd (b_n) u^\ep_t}{u^\ep_t} = \frac{1}{2} \scal{ \Opd \Big( \Big( \left\{ \lambda_n, \beta \right\} - \left\{ \beta , \lambda_n \right\} \Big) \Pi_n \Big) u^\ep_t}{u^\ep_t} \ + \ \gO{\sqrt\ep}{R}.
\]
By passing to the limit,
\begin{align*}
	\int_\R \int_{\R^2} \Xi'(t) & \tr \Big( a^W (s,0,0,0,y,\sigma,\D_y) M^t(s,\sigma) \Big) \diff\nu^t(s,\sigma) \diff t \\
	& = \int_\R \int_{\R^2} \Xi(t) \tr \Big(  \left\{ \lambda_n, \beta \right\} \Pi_n(s,0,0,0,y,\sigma,\D_y) M^t(s,\sigma) \Big) \diff\nu^t(s,\sigma) \diff t, \\
	& = \int_\R \int_{\R^2} \Xi(t) \tr \Big(  \left\{ \lambda_n, \beta \right\} \Pi_n(s,0,0,0,y,\sigma,\D_y) \Big) \diff\nu^t_n(s,\sigma) \diff t \\
	& = \int_\R \int_{\R^2} \Xi(t) \left\{ \lambda_n, \beta \right\}(s,\sigma) \diff\nu^t_n(s,\sigma) \diff t, \\
	& = - \int_\R \int_{\R^2} \Xi(t) \beta(s,\sigma) \left\{ \lambda_n, \diff\nu^t_n \right\}(s,\sigma) \diff t.
\end{align*}
Then $\nu_n^t$ satisfies the equation~\eqref{eq:propagation_measure_scale_ep} and is notably continuous on $\R$.
It remains to prove that the value at zero of $\nu_n^t$ coincides with the measure of the initial condition $\nu_{n,0}$.
Let $T \in \R$, on the one hand,
\[
    \int_0^T \frac{\mathrm d}{\mathrm{dt}} \scal{\Opd (b_n) u^\ep_t}{u^\ep_t} \diff t = \scal{\Opd (b_n) u^\ep_T}{u^\ep_T} - \scal{\Opd (b_n) u^\ep_0}{u^\ep_0}.
\]
On the other hand,
\begin{align*}
    \left| \int_0^T \frac{\mathrm d}{\mathrm{dt}} \scal{\Opd (b_n) u^\ep_t}{u^\ep_t} \diff t \right| & \leqslant \frac{1}{2} \int_0^T \left| \scal{ \Opd \Big( \left\{ \lambda_n, \beta \right\} \Pi_n \Big) u^\ep_t}{u^\ep_t} \right| \diff t \ + \ \gO{\sqrt\ep}{T,R}, \\
    & \leqslant \frac{C\prth{\left\{ \lambda_n, \beta \right\} \Pi_n}T}{2} \norme{u^\ep_0}_{\L}^2 \ + \ \gO{\sqrt\ep}{T,R}.
\end{align*}
where the last inequality holds because of the Calder\'on-Vaillancourt Theorem~\eqref{eq:Calderon_Vaillancourt} and the preservation of the $\L$-norm of the solution.
Then, by passing to the limit, we have
\begin{equation*}
\begin{split}
    \limsup_{T \to 0} \limsup_{\delta \to 0} \limsup_{R \to +\infty} \lim_{\ep \to 0} \mathcal I^\ep_{\varphi_t^\ep} \left( a_{\ep,R} \right) & = \int_{\R^2} \tr \Big( a^W (s,0,0,0,y,\sigma,\D_y) \Pi_n(s,\sigma) \Big) \diff\nu_n^0(s,\sigma), \\
    & = \limsup_{T \to 0} \limsup_{\delta \to 0} \limsup_{R \to +\infty} \lim_{\ep \to 0} \mathcal I^\ep_{\varphi_0^\ep} \left( a_{\ep,R} \right).
\end{split}
\end{equation*}
Therefore, for all $n \in \Z, \nu_n^t$ is continuous and the value at zero coincides with the initial condition.
Moreover, the equality $M^t \diff \nu^t = \sum_{n \in \Z} \Pi_n \diff \nu_n^t$ holds as trace class operator so, thanks to Lebesgue dominated convergence Theorem, $M^t\diff\nu^t$ is continuous on $\R$.
\end{proof}

\begin{remark}[Zero value of the two-scale semiclassical measure at finite distance] \label{rem:difference_zero_measure} We can construct the initial data inspired by the proof of the Lemma~\ref{lem:mesure opératorielle} where we construct semiclassical measures for families valued in a separable Hilbert space.

\bigskip

\noi Let $\prth{f^\ep}_{\ep > 0}$ be a bounded family of $\L\prth{\R_s,\C}$ with the knowledge of $\diff\rho$, the semiclassical measure of $\prth{f^\ep}_{\ep > 0}$ at the scale $\sqrt\ep$.
The function $(s,\sigma) \mapsto g_1^{s,\sigma}$ is smooth, bounded and has all derivatives bounded uniformly in $y$, therefore we can consider its Weyl $\sqrt\ep$-quantization according to the variable $s$.
Let us consider initial data $\prth{\psi_0^\ep}_{\ep > 0}$ such that for all $(s,y) \in \R \times \mathrm I_{1/2}$,
\[
    u_0^\ep(s,y) = g_0^s(y)f(s) + \Big( \mathrm{Op}_{\sqrt\ep}^{(s)} \big( g_1^{s,\sigma}(y) \big)f \Big)(s),
\]
where the Weyl $\sqrt\ep$-quantization of $g_1^{s,\sigma}$ is taken only in the $s$ variable.
Then the two-scale semiclassical measure at finite distance of this initial data is for all $(s,\sigma) \in \R^2$,
\[
    \frac{1}{2\pi} \Big( \Pi_0(s) + \Pi_1(s,\sigma) + \Pi_{1,0}(s) + \Pi_{0,1}(s,\sigma) \Big) \diff\rho(s,\sigma),
\]
where for all $(k,\ell) \in \ens{0,1}^2$ and for all $(s,\sigma) \in \R^2, \Pi_{k,\ell}(s,\sigma) = g^{s,\sigma}_k \otimes g^{s,\sigma}_\ell$.
Notably, it does not commute with $\mathrm T_\E(s,\sigma)$.
\end{remark}

\subsubsection{Two-scale Wigner measure at infinity } \label{sub:M}

We prove here the second point of Theorem~\ref{thm:main_result}.
We consider the symbol of $\mathrm H^\E_\ep$ according to the second-scale quantization~\eqref{eq:secondquantization}, given by for all $(s,y) \in \R^2$ and $z \coloneqq (z_y, z_\sigma, z_\eta) \in \R^3$,
\[
    \mathcal T_\E^\infty (s, y, z) \coloneqq z_y \int_0^1 \nabla \m \Big(\Phi(s,\tau y)\Big) \cdot \mathbf{n}(s) \diff\tau \sigma_1 + z_\eta \sigma_2 - \frac{ z_\sigma }{1 + y \kappa(s)} \sigma_3,
\]
which has two eigenvalues :
\[
    \lambda_\pm^\E(s,y,z) \coloneqq \pm \sqrt{ \frac{ 1 }{(1 + y \kappa(s))^2} z_\sigma^2 + z_y^2 \abs{\int_0^1 \nabla \m \Big(\Phi(s,\tau y)\Big) \cdot \mathbf{n}(s) \diff\tau}^2 + z_\eta^2 },
\]
which are homogeneous functions of order 1 in $z$.
The eigenprojectors $\Pi_\pm^\E(s,y,z)$ are defined by
\[
    \Pi_\pm^\E (s,y,z) \coloneqq \frac{1}{2} \mathrm{Id} + \frac{1}{2\lambda_\pm^\E(s,y,z)} \mathcal T_\E^\infty (s, y, z)
\]
and are homogeneous functions of order 0 in $z$.

\bigskip

\noi Moreover, we denote the principal symbol of the operator $\H^\E$ associated with the second-scale quantization~\eqref{eq:secondquantization} by for all $(s,y) \in \R^2$ and for all $z \coloneqq (z_y, z_\sigma, z_\eta) \in \R^3$,
\[
    \mathrm T^\infty_\E (s, z) = \mathcal T_\E^\infty (s, 0, z) = \matp{ z_\sigma & z_y r(s)^2 - iz_\eta \\ z_y r(s)^2 + iz_\eta & -z_\sigma},
\]
where $\mathrm T^\infty_\E$ was defined in~\eqref{eq:Hamiltonian_pp_scale_delta}.
Its eigenvalues are
\[
    \lambda_\pm^\infty(s,z) \coloneqq \pm \sqrt{ z_\sigma^2 + z_y^2r(s)^4 + z_\eta^2 },
\]
and the associated eigenprojectors $\Pi_\pm^\infty(s,z)$ were defined in~\eqref{eq:proj_infty}. 

\bigskip

\noi Then, the equation~\eqref{eq:Dirac_change_variable} writes on $\R_t \times \R_s \times \mathrm I_{1/2}$
\[
    \sqrt\ep \D_t \varphi^\ep_t + \Opp \Big( \mathcal T_\E^\infty \Big) \varphi^\ep_t = 0.
\]

\begin{proof}[Proof of the second point of Theorem~\ref{thm:main_result}]
Assume Assumption~\ref{ass:transversality} and Assumption~\ref{ass:tubularneighborhood}.
For convenience, we assume that the sequence $\prth{\ep}_{\ep > 0}$ realizes the measure $M^t_\infty \diff\nu^t_\infty$.
Let $a \in \mathcal{A}$ and $R, \delta$ as defined in~\eqref{eq:relation_ep_delta}.
We consider the part~\eqref{eq:observable_delta} of $a$ divided by $\lambda^\infty_+$.
The symbol $\mathcal T_\E^\infty$ is homogeneous in $z$ of order $1$ and the symbol $a_{\ep,\delta}^R / \lambda^\infty_+$ is homogeneous in $z$ of order $-1$ at infinity.
Because of~\eqref{eq:observable_delta}, we consider symbols on a subset of $\ens{z \in \R^3 | \frac{R}{2} \leqslant \abs{z} \leqslant \frac{2\delta}{\sqrt\ep} }$.
Moreover, we have the following properties.
\begin{enumerate}
    \item $\ens{\mathcal T_\E^\infty, \frac{a_{\ep,\delta}^R}{\lambda^\infty_+} }_{z_y,z_\eta}$ is homogeneous in $z$ of order $-2$.
    \item $\ens{\frac{a_{\ep,\delta}^R}{\lambda^\infty_+}, \mathcal T_\E^\infty}_{s,z_\sigma}$ is homogeneous in $z$ of order $-1$.
    \item $\ens{\mathcal T_\E^\infty, \frac{a_{\ep,\delta}^R}{\lambda^\infty_+} }_{(s,y),(\sigma,\eta)}$ is homogeneous in $z$ of order $0$.
\end{enumerate}
Then, according to standard semiclassical calculus,
\begin{align*}
    \Opp \left( \ens{\mathcal T_\E^\infty, \frac{a_{\ep,\delta}^R}{\lambda^\infty_+} }_{z_y,z_\eta} - \ens{\frac{a_{\ep,\delta}^R}{\lambda^\infty_+}, \mathcal T_\E^\infty}_{z_y,z_\eta} \right) & = \gO{\frac{1}{R^2}}{}, \\
    \sqrt\ep \Opp \left( \ens{\mathcal T_\E^\infty, \frac{a_{\ep,\delta}^R}{\lambda^\infty_+} }_{s,z_\sigma} - \ens{\frac{a_{\ep,\delta}^R}{\lambda^\infty_+}, \mathcal T_\E^\infty}_{s,z_\sigma} \right) & = \gO{\sqrt\ep}{R}, \\
    \ep \Opp \left( \ens{\mathcal T_\E^\infty, \frac{a_{\ep,\delta}^R}{\lambda^\infty_+} }_{(s,y),(\sigma,\eta)} - \ens{\frac{a_{\ep,\delta}^R}{\lambda^\infty_+}, \mathcal T_\E^\infty}_{(s,y),(\sigma,\eta)} \right) & = \gO{\ep}{R}.
\end{align*}
So
\begin{align*}
    \left[ \Opp \left( \frac{a_{\ep,\delta}^R}{\lambda^\infty_+} \right), \Opp \prth{ \mathcal T_\E^\infty } \right] = & \ \Opp \left( \crch{ \frac{a_{\ep,\delta}^R}{\lambda^\infty_+} ,\mathcal T_\E^\infty} \right) \\
    & + \frac{1}{2\i} \Opp \left( \ens{\frac{a_{\ep,\delta}^R}{\lambda^\infty_+}, \mathcal T_\E^\infty}_{z_y,z_\eta} - \ens{\mathcal T_\E^\infty, \frac{a_{\ep,\delta}^R}{\lambda^\infty_+} }_{z_y,z_\eta} \right) + \gO{\frac{1}{R^4}}{} \\
    & + \frac{\sqrt\ep}{2\i} \Opp \left( \ens{\frac{a_{\ep,\delta}^R}{\lambda^\infty_+}, \mathcal T_\E^\infty}_{s,z_\sigma} - \ens{\mathcal T_\E^\infty, \frac{a_{\ep,\delta}^R}{\lambda^\infty_+} }_{s,z_\sigma} \right) + \gO{\ep}{R}.
\end{align*}
\begin{itemize}
    \item[$\bullet$] Step 1 : Structure of $M_\infty^t$
\end{itemize}
\noi Let $\Xi \in \mathscr C^\infty_c (\R)$, then
\begin{multline*}
    \int_\R \Xi(t) \i \sqrt\ep \frac{\mathrm d}{\mathrm{dt}} \scal{ \Opp \left( \frac{a_{\ep,\delta}^R}{\lambda^\infty_+} \right) \varphi^\ep_t}{\varphi^\ep_t} \diff t = \int_\R \Xi(t) \scal{ \Opp \left( \crch{ \frac{a_{\ep,\delta}^R}{\lambda^\infty_+} ,\mathcal T_\E^\infty} \right) \varphi^\ep_t}{\varphi^\ep_t} \diff t \\
    + \frac{\sqrt\ep}{2\i} \int_\R \Xi(t) \scal{ \Opp \left( \ens{\frac{a_{\ep,\delta}^R}{\lambda^\infty_+}, \mathcal T_\E^\infty}_{s,z_\sigma} - \ens{\mathcal T_\E^\infty, \frac{a_{\ep,\delta}^R}{\lambda^\infty_+} }_{s,z_\sigma} \right) \varphi^\ep_t}{\varphi^\ep_t} \diff t + \ \gO{\ep}{R} + \gO{\frac{1}{R^2}}{} \\
    = -\i\sqrt\ep \int_\R \Xi'(t) \scal{ \Opp \left( \frac{a_{\ep,\delta}^R}{\lambda^\infty_+} \right) \varphi^\ep_t}{\varphi^\ep_t} \diff t.
\end{multline*}
But $\crch{ \frac{a_{\ep,\delta}^R}{\lambda^\infty_+}, \mathcal T_\E^\infty} = \crch{ a, \frac{\mathcal T_\E^\infty}{\lambda^\infty_+} }_{\ep,\delta}^R$ and $\crch{ a, \frac{\mathcal T_\E^\infty}{\lambda^\infty_+} } \in \mathcal{A}$, then by passing to the limit, we deduce that for all $\Xi \in \mathscr C^\infty_c (\R)$ and for all $a\in \mathcal{A}$,
\begin{align*}
    0 & = \int_{\R^2 \times \S^2} \Xi(t) \mathrm{tr} \Bigg( \left[a, \frac{\mathcal T_\E^\infty}{\lambda^\infty_+} \right]_\infty (s, 0, 0, 0, \omega) M_\infty^t (s,\omega) \Bigg) \diff \nu_\infty^t (s,\omega) \diff t \\
    & = \int_{\R^2 \times \S^2} \Xi(t) \mathrm{tr} \Bigg( \left[a_\infty(s, 0, 0, 0, \omega), \frac{1}{\lambda^\infty_+(s,\omega)} \mathcal T_\E^\infty(s,0,\omega) \right] M_\infty^t (s,\omega) \Bigg) \diff \nu_\infty^t (s,\omega) \diff t \\
    & = \int_{\R^2 \times \S^2} \Xi(t) \frac{1}{\lambda^\infty_+(s,\omega)} \mathrm{tr} \Bigg( a_\infty(s, 0, 0, 0, \omega) \Big[ \mathrm T^\infty_\E(s,\omega), M_\infty^t (s,\omega) \Big] \Bigg) \diff \nu_\infty^t (s,\omega)\diff t.
\end{align*}
So, for almost all $t\in\R$, we have $\nu_\infty^t$ a.e. $\Big[ \mathrm T^\infty_\E(s,\omega), M_\infty^t (s,\omega) \Big] = 0$ which implies~\eqref{eq:measure_scale_delta}.
\begin{itemize}
    \item[$\bullet$] Step 2 : Equation of $\diff\nu_\infty^t$.
\end{itemize}
\noi Let $\beta \in \mathcal{A}$ and $a$ defined by for all $(s,y,\sigma,\eta,z) \in \R^7$,
\[
    a(s,y,\sigma,\eta, z) \coloneqq \lambda^\infty_\pm(s,z)\beta(s,y,\sigma,\eta,z) \Pi_\pm^\E(s,y,z).
\]
The form of $a$ implies $\crch{a,\mathcal T_\E^\infty} = 0$. Then
\begin{align*}
    \sqrt\ep \frac{\diff}{\mathrm{dt}} \scal{ \Opp \left( a_{\ep,\delta}^R \right) \varphi^\ep_t}{\varphi^\ep_t} = & \ \frac{1}{2} \scal{ \Opp \left( \ens{\mathcal T_\E^\infty, a_{\ep,\delta}^R }_{z_y,z_\eta} - \ens{a_{\ep,\delta}^R, \mathcal T_\E^\infty}_{z_y,z_\eta} \right) \varphi^\ep_t}{\varphi^\ep_t} + \gO{\frac{1}{R^2}}{} \\
    & + \frac{\sqrt\ep}{2} \scal{ \Opp \left( \ens{\mathcal T_\E^\infty, a_{\ep,\delta}^R }_{s,z_\sigma} - \ens{a_{\ep,\delta}^R, \mathcal T_\E^\infty}_{s,z_\sigma} \right) \varphi^\ep_t}{\varphi^\ep_t} + \gO{\ep}{R}.
\end{align*}
Because of the homogeneity in $z$, we have $\ens{\mathcal T_\E^\infty, a}_{z_y,z_\eta} \in \mathcal{A}$. Moreover,
\begin{align*}
    \Opp \left( \ens{\mathcal T_\E^\infty, a_{\ep,\delta}^R }_{z_y,z_\eta} \right) = & \ \Opp \left( \left( \ens{\mathcal T_\E^\infty, a}_{z_y,z_\eta} \right)_{\ep,\delta}^R \right) + \frac{1}{R}\Opp \left( a \chi_{\frac{\delta}{\sqrt\ep}} \ens{\mathcal T_\E^\infty, \chi^R}_{z_y,z_\eta} \right) \\
    & + \frac{\sqrt\ep}{\delta} \Opp \left( a \chi^R \ens{\mathcal T_\E^\infty, \chi_{\frac{\delta}{\sqrt\ep}} }_{z_y,z_\eta} \right).
\end{align*}
The function $\partial\chi^R$ is compactly supported in the ring $\ens{z \in \R^3 | \frac{R}{2} \leqslant \abs{z} \leqslant R}$ and the function $a \chi_{\frac{\delta}{\sqrt\ep}} \ens{\mathcal T_\E^\infty, \chi^R}_{z_y,z_\eta}$ belongs to $\mathcal{A}$ so
\[
    \left( a \chi_{\frac{\delta}{\sqrt\ep}} \ens{\mathcal T_\E^\infty, \chi^R}_{z_y,z_\eta} \right)_\infty = 0.
\]
The operator $\Opp \left( a \chi^R \ens{\mathcal T_\E^\infty, \chi_{\frac{\delta}{\sqrt\ep}} }_{z_y,z_\eta} \right)$ is uniformly bounded in $\ep$ then
\[
     \lim_{\ep \to 0}  \frac{\sqrt\ep}{\delta} \scal{ \Opp \left( a \chi^R \ens{\mathcal T_\E^\infty, \chi_{\frac{\delta}{\sqrt\ep}} }_{z_y,z_\eta} \right) \varphi_t^\ep}{\varphi_t^\ep} = 0.
\]
As in the computation of~\eqref{eq:decomposition_proj_n}, we consider the quantity $\Pi_+^\E \ens{a, \mathcal T_\E^\infty }_{z_y , z_\eta} \Pi_+^\E$ then
\[
    \Pi_+^\E \ens{a, \mathcal T_\E^\infty }_{z_y , z_\eta} \Pi_+^\E = \lambda_+^\infty \ens{ \beta , \lambda_+^\infty }_{z_y , z_\eta} \Pi_+^\E.
\]
Then, by passing to the limit, we have
\begin{align*}
    0 & = \int_{\R^2 \times \S^2} \Xi(t) \mathrm{tr} \Bigg( \left( \lambda_+^\infty \ens{ \beta , \lambda_+^\infty }_{z_y , z_\eta} \Pi_+^\E \right)_\infty (s, 0, 0, 0, \omega) M_\infty^t (s,\omega) \Bigg) \diff \nu_\infty^t (s,\omega) \diff t, \\
    & = \int_{\R^2 \times \S^2} \Xi(t) \mathrm{tr} \Bigg( \lambda_+^\infty \ens{ \beta_\infty , \lambda_+^\infty }_{\omega_y , \omega_\eta} (s,0,0,0,\omega) \Pi_+^\infty (s,\omega) \Bigg) \diff \nu_+^t (s,\omega) \diff t, \\
    & = \int_{\R^2 \times \S^2} \Xi(t) \lambda_+^\infty \ens{ \beta_\infty , \lambda_+^\infty }_{\omega_y , \omega_\eta} (s,0,0,0,\omega) \diff \nu_+^t (s,\omega) \diff t,
\end{align*}
but
\[
    \lambda_+^\infty(s,\omega)\frac{\partial\lambda_+^\infty}{\partial z_y} (s,\omega) = \omega_yr(s)^4, \quad  \lambda_+^\infty(s,\omega)\frac{\partial\lambda_+^\infty}{\partial z_\eta} (s,\omega) = \omega_\eta, \quad (s,\omega) \in \R \times \S^2.
\]
Then
\[
    \int_{\R^2 \times \S^2} \Xi(t) \Big( \vec{V}_+(s,\omega) \cdot \nabla_{\omega} \beta_\infty (s,\omega) \Big) \diff \nu_+^t (s,\omega) \diff t = 0,
\]
with
\[
    \vec{V}_+(s,\omega) = \matp{-\omega_\eta \\ 0 \\ \omega_y r(s)^4}
\]
and $\beta_\infty$ a function on the sphere with $\nabla_\omega \beta_\infty$ is a vector field over the sphere. Therefore we can rewrite 
$\vec V_\pm \cdot \nabla_\omega \beta =\vec V_\pm^\infty \cdot \nabla_\omega \beta$ where $V_\pm^\infty$ is the component of $\vec V_\pm$ which is tangential to the sphere. The vector $\vec V_\pm^\infty$ defines a vector field over the sphere and is given by for all $(s,\omega) \in \R \times \S^2$,
\[
    \vec{V}^\infty (s,\omega) \coloneqq \vec{V}_+ (s,\omega) - \Big( \vec{V}_+ (s,\omega)\cdot \omega \Big) \omega = \matp{ \prth{ \prth{1 - r(s)^4}\omega_y^2 - 1} \omega_\eta \\ \prth{1 - r(s)^4} \omega_y\omega_\eta \omega_\sigma \\ \prth{ \prth{1-r(s)^4} \omega_\eta^2 + r(s)^4} \omega_y }.
\]
Because of the equality $\lambda^\infty_+ = - \lambda^\infty_-$, we have the same invariance for $\diff \nu_-^t$. So, $t$ almost everywhere,
\[
    \int_{\R \times \S^2} \Big( \vec{V}^\infty(s,\omega) \cdot \nabla_{\omega} \beta_\infty (s,\omega) \Big) \diff \nu_\infty^t (s,\omega) = 0.
\]
So, for almost all $(t,s,\omega) \in \R^2\times\S^2$, we have $\mathrm{div}_\omega \Big( \vec{V}^\infty(s,\omega) \diff\nu^t (s,\omega) \Big) = 0$.
\end{proof}

\subsection{The Wigner measure} \label{sec:proof_main}

\begin{proof}[Proof of Theorem \ref{thm:main_result}] Assume Assumption~\ref{ass:transversality} and Assumption~\ref{ass:tubularneighborhood}.
Let $\prth{\psi^\ep_t}_{\ep > 0}$ a uniformly bounded family in $\mathrm L^\infty \prth{ \R, \L \prth{ \R^2,\C^2}}$ solution to equation~\eqref{eq:Dirac} with normalized initial condition $\prth{\psi^\ep_0}_{\ep > 0}$ in $\L\prth{\R^2,\C^2}$.

\medskip

\noi Let us consider the family $\prth{\varphi^\ep_t}_{\ep > 0}$ defined in~\eqref{def:varphi} which is solution to equation~\eqref{eq:Dirac_change_variable}, according to Lemma~\ref{lem:normal form} and uniformly bounded in $\mathrm L^\infty \prth{ \R, \L \prth{ \R^2,\C^2}}$.
Lemma~\ref{lem: normal form measure} gives the link between any semiclassical measure of $\prth{\psi^\ep_t}_{\ep > 0}$ over the curve $\mathcal C$ and the associated measure of $\prth{\varphi^\ep_t}_{\ep > 0}$ on $\R_s$.
According to Theorem~\ref{thm:existence_two_scale_measure}, we can consider $\displaystyle{ \prth{M^t_\infty \diff\nu^t_\infty , M^t \diff\nu^t} }$ a pair of two-scale semiclassical measure of $\prth{\varphi^\ep_t}_{\ep > 0}$.
Finally, according to proofs of Section~\ref{sub:Gamma} and Section~\ref{sub:M}, we conclude the evolution of the semiclassical measure over the curve $\mathcal C$.
\end{proof}

\subsection{Proof of applications} \label{sec:proof_applications}

\begin{proof}[Proof of Corollary~\ref{cor:conservation}] Assume Assumption~\ref{ass:transversality} and Assumption~\ref{ass:tubularneighborhood}.
With the notation of Theorem~\ref{thm:main_result}, the evolution of $(\diff\mu_t,M^t\diff\nu^t)$ implies that for all $t \in \R$,
\begin{align*}
	\int_{\R^4 \backslash \mathcal{C}} \mathrm{tr} \Big( \diff\mu_t(x,\xi) \Big) & = \int_{\R^4 \backslash \mathcal{C}} \mathrm{tr} \Big( \diff\mu_0(x,\xi) \Big), \\
	\int_{\R^2} \tr \Big( M^t(s,\sigma) \Big) \diff\nu^t(s,\sigma) & = \int_{\R^2} \tr \Big( M^0(s,\sigma) \Big) \diff\nu^0(s,\sigma).
\end{align*}
By taking $\prth{\psi^\ep_0}_{\ep > 0} \in \L(\R^2,\C^2)$ normalized such that
\[
    \int_{\R^4 \backslash \mathcal{C}} \mathrm{tr} \Big( \diff\mu_0(x,\xi) \Big) + \int_{\R^2} \tr \Big( M^0(s,\sigma) \Big) \diff\nu^0(s,\sigma) = 1 = \norme{\psi^\ep_0}_2^2,
\]
we have the following equality because of the preservation of the $\L$-norm
\begin{align*}
    \norme{\psi^\ep_t}_2^2 = & \ \int_{\R^4 \backslash \mathcal{C}} \mathrm{tr} \Big( \diff\mu_t(x,\xi) \Big) + \int_{\R^2} \tr \Big( M^t(s,\sigma) \Big) \diff\nu^t(s,\sigma) + \int_{\R \times \S^2} \mathrm{tr} \Big(M^t_\infty (s,\omega) \Big) \diff\nu^t_\infty (s,\omega), \\
    = & \ 1.
\end{align*}
Therefore, for almost all $t$,
\[
    M^t_\infty (s,\omega) \diff\nu^t_\infty (s,\omega) = 0.
\]
\end{proof}

\begin{proof}[Proof of Corollary~\ref{cor:goal_reached}]
Assume Assumption~\ref{ass:transversality} and Assumption~\ref{ass:tubularneighborhood}.
According to Lemma~\ref{lem:GMMP}, we know the evolution of any semiclassical measure of $\prth{\psi^\ep_t}_{\ep > 0}$ outside the curve $\mathcal C$.
According to Theorem~\ref{thm:main_result}, by applying Remark~\ref{rem:measure_observable_phase_space}, we can describe the semiclassical measure above the interface so we conclude.
\end{proof}

\begin{proof}[Proof of Corollary~\ref{cor:linkDrouot}] Assume Assumption~\ref{ass:transversality} and Assumption~\ref{ass:tubularneighborhood}.
Let $\vec f \in \L\prth{\R^2,\C^2}$. Let $\displaystyle{ \prth{\psi_t^\ep}_{\ep > 0} }$ be the solution to equation~\eqref{eq:Dirac} with, for all $\ep \in (0,1]$ and for all $x \in \R^2$,
\[
    \psi_0^\ep(x) = \frac{1}{\sqrt\ep} \vec f \prth{ \frac{x-x_0}{\sqrt\ep} },
\]
then, applying Corollary~\ref{cor:goal_reached} and Corollary~\ref{cor:condition_M_infty_vanishing}, there exist a vanishing sequence of positive numbers $\displaystyle{(\ep_k)_{k \in \N}}$ and a two-scale semiclassical measure at finite distance $\displaystyle{M^t \diff\nu^t }$, such that for all observables $(a,\Xi) \in \mathscr C_c^\infty\prth{\R^4,\C^{2,2}} \times \mathscr C_c^\infty(\R,\C)$,
\[
    \mathcal I^{\ep_k}_\psi(\Xi,a) \underset{k \to + \infty}{\longrightarrow} \frac{1}{2} \sum_{n \in \Z} \int_{\R^3} \Xi(t) \mathrm{tr} \left( a \big( \mathbf{t}(s),0 \big) \matp{ 1 & \frac{-\sigma}{\lambda_n(s,\sigma)} e^{ -\i\theta(s)} \\ \frac{-\sigma}{\lambda_n(s,\sigma)}e^{ \i\theta(s)} & 1} \diff\nu^t_n(s,\sigma) \right) \diff t,
\]
with, for all $n\in\Z$, $\nu_n^t$ satisfies the following equation
\[
    \left\{
    \begin{array}{rll}
    \partial_t \ \nu_n^t & = & \ens{\nu_n^t, \lambda_n} \\
    \nu_n^t |_{t=0} & = & \nu_{n,0} 
    \end{array}
    \right. ,
\]
with $\nu_{n,0}= \mathrm{tr} \Big( \Pi_nM^0 \Big) \nu^0$ and $M^0 \diff\nu^0 $ is a two-scale semiclassical measure associated with the initial data.

\medskip

\noi We conclude thanks to the conservation of $\ep$-oscillation (Lemma~\ref{eposc}).
\end{proof}

\begin{proof}[Proof of Corollary~\ref{cor:below_estimate}] Assume Assumption~\ref{ass:transversality} and Assumption~\ref{ass:tubularneighborhood}. Let $\displaystyle{ \prth{\psi_0^\ep}_{\ep > 0} }$ a uniformly bounded family of $\L\prth{\R^2,\C^2}$, $\nu$ a non-negative scalar Radon measure on $\R_t \times \R_x^2$ and a vanishing sequence of positive numbers $\displaystyle{(\ep_k)_{k \in \N}}$ such that for all $(a,\Xi) \in \mathscr C_c^\infty \prth{ \Omega ,\C} \times \mathscr C_c^\infty( \R ,\C)$,
\[
    \int_\R \int_\Omega \Xi(t) a(x) \abs{\psi_t^{\ep_k}(x)}^2 \diff x \diff t \underset{k \to +\infty}{\longrightarrow} \int_\R \int_\Omega \Xi(t) a(x) \diff\nu(t,x),
\]
where $\displaystyle{ \prth{\psi_t^\ep}_{\ep > 0} }$ solves~\eqref{eq:Dirac} with $\displaystyle{ \prth{\psi_0^\ep}_{\ep > 0} }$ as initial condition.
A consequence of Theorem~\ref{thm:main_result} and Remark~\ref{rem:measure_observable_phase_space} is
\[
    \diff\nu \mathds{1}_\mathcal C \geqslant \sum_{n \in \Z} \Phi^* \left( \tr \Big( \Pi_n(s,\sigma) \Big) \diff\nu^t_n \right) \mathds{1}_\mathcal C \diff t
\]
where $\prth{\nu_n^t}_{n \in \Z}$ is defined by~\eqref{eq:propagation_measure_scale_ep}.
We conclude with $\tr \left( \Pi_n(s,\sigma) \right) = 1$.
\end{proof}

\appendix

\section{Semiclassical analysis} \label{app:A}

Our problem requires an analysis in phase space.
We consider the {\it Wigner transform} of the solution $\prth{\psi_t^\ep}_{\ep > 0}$ of equation~\eqref{eq:Dirac} and its weak limits in the space of distributions that are called {\it Wigner measures}.
We apply the main properties of the Wigner transform (introduced in Section~\ref{subsec:wigner}) to the solution of our problem in Section~\ref{app:Dirac}.

\bigskip

First, let us briefly introduce an important example and some definitions.

\begin{example} \label{example1} Let us consider $\alpha \in (0,1]$, $f \in \L\prth{\R^2,\C^2}$ and for all $x \in \R^2$,
\[
    f^\ep(x) \coloneqq \frac{1}{\ep^{\alpha/2}} f\prth{ \frac{x}{\ep^\alpha} },
\]
so $\displaystyle{\prth{f^\ep}_{\ep > 0}}$ has only one Wigner measure that depends of $\alpha$.

\bigskip

\noi If $\alpha = 1/2$, then for all $(x,\xi) \in \R^4$, $\displaystyle{\diff\mu(x,\xi) = \prth{\int_{\R^2} f(x) \otimes \overline{f(x)} \diff x } \delta_0(x) \otimes \delta_0 (\xi) }$.

\noi If $\alpha = 1$, then for all $(x,\xi) \in \R^4$, $\displaystyle{\diff\mu(x,\xi) = \widehat{f}(\xi) \otimes \overline{\widehat{f}(\xi)} \delta_0(x) \frac{\diff\xi}{2\pi} }$.
\end{example}

\noi We say a bounded family $\prth{f^\ep}_{\ep > 0}$ of $\L\prth{\R^2,\C^2}$ is $\ep$-oscillating if
\[
    \limsup_{\ep \to 0} \int_{\abs{\xi} > R/\ep} \abs{\widehat{f^\ep}(\xi)}^2_{\C^2} \diff\xi \underset{R \to +\infty}{\longrightarrow} 0.
\]
We say a bounded family $\prth{f^\ep}_{\ep > 0}$ of $\L\prth{\R^2,\C^2}$ is compact at infinity if
\[
    \limsup_{\ep \to 0} \int_{\abs{x} > R} \abs{f^\ep(x)}^2_{\C^2} \diff x \underset{R \to +\infty}{\longrightarrow} 0.
\]
If $\prth{f^\ep}_{\ep > 0}$ is $\ep$-oscillating and $\prth{\abs{f^\ep}^2_{\C^2} \diff x}_{\ep > 0}$ converges weakly to $\nu$, then for all semiclassical measure $\diff\mu$ of $\prth{f^\ep}_{\ep > 0}$,
\[
    \nu(x) = \mathrm{tr} \int_{\R^2} \diff\mu(x,\diff\xi).
\]
If $\prth{f^\ep}_{\ep > 0}$ is compact at infinity and $\prth{\abs{\widehat{f^\ep}}^2_{\C^2} \diff \xi}_{\ep > 0}$ converges weakly to $\widetilde{\nu}$, then for all semiclassical measure $\diff\mu$,
\[
    \widetilde{\nu}(\xi) = \mathrm{tr}\int_{\R^2} \diff\mu(\diff x,\xi).
\]
Now, we apply these first properties on solutions to equation~\eqref{eq:Dirac}.

\subsection{Application to Dirac equation} \label{app:Dirac}

In this section, we consider $\prth{\psi^\ep_t}_{\ep > 0}$ the solution to equation~\eqref{eq:Dirac} with the normalized initial condition $\prth{\psi^\ep_0}_{\ep > 0}$. We revisit in this paragraph the notions introduced in the preceding one for this special family.

\begin{lemma}[Conservation of $\ep$-oscillating] \label{eposc} Let $T>0$.  If $\prth{\psi^\ep_0}_{\ep > 0}$ is $\ep$-oscillating then $\prth{\psi^\ep_t}_{\ep > 0}$ is  $\ep$-oscillating uniformly for $t \in [-T,T]$. In other words, uniformly for $t \in [-T,T]$,
\[
    \limsup_{\ep \to 0} \int_{\abs{\xi} > R/\ep} \abs{\widehat{\psi^\ep_t}(\xi)}^2 \diff\xi \underset{R \to +\infty}{\longrightarrow} 0.
\]
\end{lemma}

\begin{proof}[Proof of Lemma~\ref{eposc}] Let $\prth{\psi^\ep_0}_{\ep > 0}$ be $\ep$-oscillating. The proof relies on the analysis of
\[
    \psi^{\ep,R}_t \coloneqq \chi^R(\ep \D_x) \psi_t^\ep,
\]
where $\displaystyle{\chi^R(\xi) = \chi \prth{\frac{\xi}{R}}}$ with $\chi \in \mathscr C ^\infty\prth{\R^2,[0,1]}, \chi \equiv 0$ for $\abs{\xi}< 1$ and $\chi \equiv 1$ for $\abs{\xi} > 2$. We observe that $\psi^{\ep,R}_t$ solves
\[
    \i\ep\partial_t\psi^{\ep,R}_t = \H \psi^{\ep,R}_t + \crch{\chi_R(\ep \D_x), \matp{\m(x)&0\\0&-\m(x)}} \psi^{\ep,R}_t.
\]
Nevertheless, the derivatives of $\m$ are bounded therefore there exists $C>0$ such that
\[
    \norme{\frac{1}{\ep} \crch{\chi^R(\ep \D_x), \matp{\m(x)&0\\0&-\m(x)}}}_{\mathcal L\prth{ \L\prth{\R^2,\C^2} }} \leqslant \frac{C}{R}.
\]
Then, by an energy argument, for all $t \in [-T,T]$,
\[
    \frac{\diff}{\diff \mathrm t} \norme{\psi^{\ep,R}_t}_{\L\prth{\R^2,\C^2}} \leqslant \frac{C}{R},
\]
so,
\begin{align*}
    \limsup_{\ep\rightarrow 0}\norme{\psi^{\ep,R}_t}_{\L\prth{\R^2,\C^2}} & \leqslant \limsup_{\ep\rightarrow 0}\norme{\psi^{\ep,R}_0}_{\L\prth{\R^2,\C^2}} + \frac{C}{R} \abs{T}, \\
    & \underset{R \to +\infty}{\longrightarrow} 0.
\end{align*}
because $\prth{\psi^\ep_0}_{\ep > 0}$ is uniformly bounded in $\L\left(\R^2,\C^2\right)$ and $\ep$-oscillating.
\end{proof}

The $\ep$-oscillation conservation property justifies the semiclassical approach to solutions to equation~\eqref{eq:Dirac} for analyzing the limits of the quantities~\eqref{eq:Wignerdistribution}.
Note that similar arguments show that there is also no loss of mass at infinity in configuration space in the sense that if $\prth{\psi_0^\ep}_{\ep > 0}$ is compact at infinity, then $\prth{\psi^\ep_t}_{\ep > 0}$ is compact at infinity uniformly in time $t \in [-T,T]$, for some $T>0$.

\begin{lemma}[Conservation of compactness at infinity] \label{compactatinf} Let $T>0$.  If $\prth{\psi^\ep_0}_{\ep > 0}$ is compact at infinity then $\prth{\psi^\ep_t}_{\ep > 0}$ solution to~\eqref{eq:Dirac} is compact at infinity uniformly for $t\in [-T,T]$.
\end{lemma}

\begin{proof}[Proof of Lemma~\ref{compactatinf}] Let $\prth{\psi^\ep_0}_{\ep > 0}$ be compact at infinity. The proof relies on the analysis of
\[
    \psi_R^\ep(t) \coloneqq \chi_R(x) \psi_t^\ep,
\]
where $\displaystyle{\chi_R(x) = \chi \prth{\frac{x}{R}}}$ with $\chi \in \mathscr C ^\infty\prth{\R^2,[0,1]}, \chi \equiv 0$ for $\abs{x}< 1$ and $\chi \equiv 1$ for $\abs{x} > 2$. We observe that $\psi_R^\ep$ solves
\[
    \i\ep\partial_t\psi_R^\ep = \H \psi_R^\ep + \ep \crch{\chi_R(x), \matp{0& \D_1 - \i\D_2 \\ \D_1 + \i\D_2 & 0} } \psi_t^\ep.
\]
Nevertheless, the first derivatives of $\chi$ are bounded therefore there exists $C>0$ such that
\[
    \norme{\crch{\chi_R(x), \matp{0& \D_1 - \i\D_2 \\ \D_1 + \i\D_2 & 0}}}_{\mathcal L(\L\prth{\R^2,\C^2})} \leqslant \frac{C}{R}.
\]
Then, by an energy argument,
\[
    \frac{\diff}{\diff \mathrm t} \norme{\psi_R^\ep(t)}_{\L\prth{\R^2,\C^2}} \leqslant \frac{C}{R} \norme{\psi^\ep_0}_{\L\prth{\R^2,\C^2}},
\]
so
\begin{align*}
    \norme{\psi_R^\ep(t)}_{\L\prth{\R^2,\C^2}} & \leqslant \norme{\psi_R^\ep(0)}_{\L\prth{\R^2,\C^2}} + \frac{C}{R} \abs{t} \norme{\psi^\ep_0}_{\L\prth{\R^2,\C^2}}, \\
    & \underset{R \to +\infty}{\longrightarrow} 0.
\end{align*}
\end{proof}

As a consequence, if $\prth{\psi_0^\ep}_{\ep>0}$ is normalized, $\ep$-oscillating and compact at infinity, then for all $\chi \in \mathscr C_c^\infty\prth{\R^2,[0,1]}$, with $\chi \equiv 1$ for $\abs{x}< 1$, $\chi \equiv 0$ for $\abs{x} > 2$,
\[
    \lim_{\ep \to 0} \scal{\Op \prth{\chi \prth{\frac{x}{R}} \chi \prth{\frac{\ep \D_x}{R}} } \psi^\ep_0}{\psi^\ep_0}_{\L\prth{\R^2,\C^2}} \; \underset{R \to +\infty}{\longrightarrow} \; \lim_{\ep \to 0} \norme{\psi_0^\ep}^2_{\L\prth{R^2,\C^2}} = 1.
\]
So for all semiclassical measure $\diff\mu$ of $\prth{\psi_0^\ep}_{\ep>0}$ normalized, $\ep$-oscillating and compact at infinity,
\[
    \mathrm{tr} \int_{\R^4} \diff\mu(x, \xi) = 1.
\]

\subsection{Semiclassical measures for families valued in a separable Hilbert space} \label{app:separable}

We will use the more general framework of~\cite{Gerardmicrolocal} and consider families valued in a separable Hilbert space $\mathcal H$.
More precisely, in this section, $\mathcal H = \L \prth{\R,\C^2}$.
Note that in Section~\ref{subsec:wigner}, $\mathcal H = \C^2$.

\begin{lemma}[Operator valued measure] \label{lem:mesure opératorielle}
Let $\displaystyle{\prth{f^\ep}_{\ep > 0}}$ a bounded family in $\L \prth{\R,\mathcal H}$, then there exists a vanishing sequence of positive numbers $\displaystyle{(\ep_k)_{k > 0}}$ and a positive operator valued measure $M\diff\nu$ such that for all $a \in \mathscr C_c^\infty \prth{ \R^2, \mathcal{K} \prth{\mathcal H} }$,
\[
    \scal{ \mathrm{Op}_{\ep_k}(a) f^{\ep_k} }{ f^{\ep_k}} \underset{k \to +\infty}{\longrightarrow} \int_{\R^2} \tr_{\mathcal L(\mathcal H)} \Big( a(x,\xi) M(x,\xi) \Big) \diff\nu(x,\xi),
\]
with $\mathcal{K} \prth{\mathcal H}$ the set of compact operators on $\mathcal H$.
\end{lemma}

\begin{proof}[Proof of Lemma~\ref{lem:mesure opératorielle}] Let $\mathcal H$ be a separable Hilbert space, $\displaystyle{\prth{f^\ep}_{\ep > 0}}$ uniformly bounded in $\L \prth{\R,\mathcal H}$, $(h_n)_{n \in \N}$ a Hermitian basis of $\mathcal H$ and $\prth{\prth{f^\ep_n}_{\ep > 0}}_{n \in \N}$ families of $\L\prth{\R}$ such that $x$ a.e.
\[
    f^\ep(x) = \sum_{n \in \N} \underbrace{f^\ep_n(x)}_{\in \C} h_n.
\]
Then the families $\prth{f^\ep_n}_{\ep > 0}$ satisfy, for some constant $C > 0$,
\begin{equation} \label{eq:estimate_family}
    \sum_{n \in \N} \norme{f^\ep_n}_{\L(\R)}^2 = \norme{f^\ep}_{\L\prth{\R,\mathcal H}}^2 \leqslant C.
\end{equation}
For all $n \in \N$, $\prth{f^\ep_n}_{\ep>0}$ are bounded in $\L\prth{\R}$. By a diagonal extraction process, we can find a sequence $\displaystyle{(\ep_k)_{k > 0}}$ and a family of Radon measures $\prth{ \diff\mu_{n,m} }_{(n,m) \in \N^2}$ such that for all $b \in \mathscr C_c^\infty \prth{\R^2}$ and for all $(n,m) \in \N^2$,
\[
    \scal{ \mathrm{Op}_{\ep_k}(b) f^{\ep_k}_n }{ f^{\ep_k}_m} \underset{k \to +\infty}{\longrightarrow} \int_{\R^2} \mathrm{tr}  \Big( b(x,\xi) \diff\mu_{n,m}(x,\xi) \Big),
\]
with $\displaystyle{ \diff\mu_{n,n} \prth{ \R^2 } \leqslant \limsup_{k \to +\infty} \norme{f^{\ep_k}_n}_{\L(\R)}^2 }$ and $\diff\mu_{n,m}$ absolutely continuous with respect to $\diff\mu_{n,n}$ and $\diff\mu_{m,m}$. Moreover, for all $(n,m) \in \N^2$, $\diff\mu_{n,n}$ is a positive Radon measure and $\diff\mu_{m,n} = \overline{\diff\mu_{n,m}}$.

\bigskip

\noi First of all, we construct the operator valued measure $M\diff\nu$.
By~\eqref{eq:estimate_family},
\[
    \limsup_{k \to +\infty} \sum_{n \in \N} \norme{f^{\ep_k}_n}_{\L}^2 < +\infty,
\]
then $\displaystyle{ \nu \coloneqq \sum_{n \in \N} \mu_{n,n} }$ defines a positive Radon measure of finite mass.
By the Radon-Nykodym Theorem, there exists a sequence of $\diff\nu$-measurable functions $\prth{M_{n,m}}_{(n,m) \in \N^2}$ such that $\forall (n,m) \in \N^2$,
\[
    \diff\mu_{n,m} = M_{n,m} \diff\nu
\]
with $0 \leqslant M_{n,n} \leqslant 1$ and $M_{n,m} = \overline{M_{m,n}}$ $\diff\nu$-almost everywhere. Let us define the operator valued function $M : (x,\xi) \mapsto M(x,\xi)$ such that, for all $(n,m) \in \N^2$, for all $(x,\xi) \in \R^2$,
\[
    \scal{M(x,\xi)h_n}{h_m}_\mathcal H =M_{n,m}(x,\xi).
\]
Because of the positivity of $\prth{\diff\mu_{n,n}}_{n \in \N}$ and the adjoint relation of $\prth{\diff\mu_{n,m}}_{(n,m) \in \N^2}$, $M (x,\xi)$ is a positive self-adjoint operator and its trace satisfies
\[
    \int_{\R^2} \tr_{\mathcal L(\mathcal H)} \Big( M(x,\xi) \Big) \diff\nu(x,\xi) = \sum_{n \in \N} \int_{\R^2} M_{n,n}(x,\xi) \diff\nu(x,\xi) = \sum_{n \in \N} \diff\nu_{n,n} \prth{\R^2} < +\infty.
\]
Therefore, the pair $\prth{M,\diff\nu}$ generates an operator measure $M\diff\nu$. We now have to prove that $M\diff\nu$ plays the expected role.

\bigskip

\noi Let $a \in \mathscr C_c^\infty \prth{ \R^2, \mathcal{K} \prth{\mathcal H} }$, $\Pi_{\leqslant n}$ the projector on the finite dimensional subspace of $\mathcal H$ generated by $(h_m)_{m \leqslant n}$ and consider the finite rank operator-valued function $a_n \coloneqq \Pi_{\leqslant n} a \Pi_{\leqslant n}$ which belongs to $\mathscr C_c^\infty \prth{ \R^2, \mathcal{K} \prth{\mathcal H}}$, then
\begin{equation} \label{eq:observable_convergence}
    \sup_{(x,\xi) \in \supp(a)} \norme{a(x,\xi) - a_n(x,\xi)}_{\mathcal L \prth{\mathcal H}} \underset{n \to +\infty}{\longrightarrow} 0.
\end{equation}
By passing to the limit, thanks to the construction of $\prth{\ep_k}_{k \in \N}$,
\begin{align*}
	\scal{ \mathrm{Op}_{\ep_k}\prth{a_n} f^{\ep_k} }{ f^{\ep_k}} & = \sum_{m,\ell \leqslant n} \scal{ \mathrm{Op}_{\ep_k}\Big( \scal{a_n h_m}{h_\ell}_\mathcal H \Big) f^{\ep_k}_m }{ f^{\ep_k}_\ell}_{\L\prth{\R,\C^2}}  \\
	& \underset{k \to +\infty}{\longrightarrow} \sum_{m,\ell \leqslant n} \int_{\R^2} \scal{a_n(x,\xi) h_m}{h_\ell}_\mathcal H \diff\nu_{m,\ell}(x,\xi) \\
	& = \sum_{m,\ell \leqslant n} \int_{\R^2} \scal{a_n h_m}{h_\ell}_\mathcal H M_{m,\ell} \diff\nu(x,\xi) \\
	& = \int_{\R^2} \tr_{\mathcal L(\mathcal H)} \Big( a_n(x,\xi) M(x,\xi) \Big) \diff\nu(x,\xi).
\end{align*}
Then, with~\eqref{eq:observable_convergence} and Lebesgue dominated convergence Theorem, we have
\[
    \int_{\R^2} \tr_{\mathcal L(\mathcal H)} \Big( \big( a(x,\xi) - a_n(x,\xi) \big) M(x,\xi) \Big) \diff\nu(x,\xi) \underset{n \to +\infty}{\longrightarrow} 0.
\]
\end{proof}

\section{Spectral Analysis of the principal normal form operator} \label{app:spectral_theory}

In this appendix, we study the operator $\mathrm T_\E(s,\sigma)$ introduced in~\eqref{eq:Hamiltonian_pp_scale_ep}.
First, we rescale the operator by $r(s)$.
This is possible because $r$ is smooth ($\m$ is in $\mathscr C ^\infty\prth{\R^2,\R}$) and strictly non-negative (Assumption~\ref{ass:transversality}).

\bigskip

Let us define
\[
    \mathscr T \coloneqq \begin{pmatrix} \zeta & \mathfrak{a}^\dagger \\ \mathfrak{a} & -\zeta \end{pmatrix},
\]
where $\mathfrak{a} \coloneqq y + \i\D_y$ denotes the usual annihilation operator on the variable $y$.

\begin{remark} \label{rem:spectrum}
By denoting by $\zeta$ the ratio $\sigma / r(s)$, we have the following, for all $f \in \L\prth{\R,\C^2}$,
\begin{align*}
    \mathrm T_\E (s,\sigma) \Bigg[f \Big(r(s)y\Big) \Bigg] & = \crch{ \matp{ \sigma & r(s) (y - \i\D_y) \\ r(s) (y + \i\D_y) & -\sigma } f} \Big(r(s)y\Big) \\
    & = r(s) \crch{ \matp{ \zeta & y - \i\D_y \\ y + \i\D_y & -\zeta } f} \Big(r(s)y\Big).
\end{align*}
Then the link between spectra is
\[
    \Sp \Big( \mathrm T_\E (s,\sigma) \Big) = r(s) \Sp(\mathscr T).
\]
So we can study $\mathscr T$ instead of $\mathrm T_\E (s,\sigma)$.
\end{remark}

Let us recall usual properties of creation and annihilation operators.

\begin{proposition} \label{pro:creation_annihilation} We have the following.
\begin{itemize}
    \item $\mathfrak{a}\mathfrak{a}^\dagger = -\partial_x^2+x^2+1$. 
    \item $\operatorname{sp}\left(\mathfrak{a}\mathfrak{a}^\dagger\right) = \lbrace2n \mid n\in\N_{>0}\rbrace$. 
    \item The Hermite functions $(\mathfrak{h}_n)_{n\in\N}$, defined in~\eqref{eq:Hermite_function}, form an Hilbertian basis of $\L(\R,\C)$, of eigenfunctions of $\mathfrak{a}\mathfrak{a}^\dagger$, and we have 
    \[
        \mathfrak{a}\mathfrak{a}^\dagger\mathfrak{h}_n = 2(n+1) \mathfrak{h}_n, \quad n\in\N.
    \]
    \item Moreover $\mathfrak{a}$ and $\mathfrak{a}^\dagger$ satisfies the following identities
    \begin{equation*}
        \mathfrak{a}\mathfrak{h}_{n}=\sqrt{2n}\mathfrak{h}_{n-1},\quad \mathfrak{a}^\dagger \mathfrak{h}_n=\sqrt{2(n+1)}\mathfrak{h}_{n+1}.
    \end{equation*}
\end{itemize}
\end{proposition}

\noindent The next Proposition consists to diagonalize, at fixed $(s,\sigma)$, the operator $\mathrm T_\E (s,\sigma)$ on $\L(\R_y,\C^2)$ since $(\mathrm T_\E (s,\sigma),\mathcal{B}^1(\R_y,\C^2))$ is self-adjoint with compact resolvent.

\begin{theorem}[Spectrum of $\mathrm T_\E(s,\sigma)$] The spectrum of $\mathrm T_\E (s,\sigma)$ consists in multiplicity one eigenvalues
\[
    \Sp \Big( \mathrm T_\E (s,\sigma) \Big) = \left\{ \lambda_0(\sigma) \right\} \bigcup \left\{ \lambda_n(s,\sigma) \ | \ n \in \Z_{\neq 0} \right\}, \quad (s,\sigma) \in \R^2,
\]
with, for all $n \in \Z_{\neq 0}$,
\[
    \lambda_0(\sigma) \coloneqq \sigma, \quad \lambda_n(s,\sigma) \coloneqq \sgn(n) \sqrt{\sigma^2 + 2|n|r(s)^2}.
\]
Moreover, an Hilbertian basis of $\L \left( \R_y, \C^2 \right) $, composed of eigenfunctions of $\mathrm T_\E (s,\sigma)$, is given by $\displaystyle{\left(g_n^{s,\sigma}\right)_{n \in \Z}}$
\begin{align*}
    g_0^s (y) & \coloneqq \frac{1}{\sqrt{r(s)}}\begin{pmatrix} \mathfrak{h}_0\big(r(s)y\big) \\ 0 \end{pmatrix}, \\
    g_n^{s,\sigma} (y) & \coloneqq \alpha_n(s,\sigma) \begin{pmatrix} \frac{\sqrt{2|n|} r(s)}{\lambda_n(\sigma) - \sigma} \mathfrak{h}_{|n|} \big(r(s)y\big) \\ \mathfrak{h}_{|n|-1}\big(r(s)y\big) \end{pmatrix}, \quad \quad n \in \Z_{\neq 0} \\
    \mbox{with} \quad \quad \alpha_n (s,\sigma) & \coloneqq \frac{1}{\sqrt{2r(s)}} \sqrt{ 1 - \frac{\sigma}{\lambda_n(s,\sigma)} },
\end{align*}
and $\left( \mathfrak{h}_n \right)_{n \in \N}$ are Hermite functions defined in~\eqref{eq:Hermite_function}.
Moreover, the projectors over the eigenspaces of $\mathrm T_\E (s,\sigma)$ defined by
\[
    \Pi_0 (s,\sigma) \coloneqq g_0^s \otimes g_0^s, \quad \quad \Pi_n (s,\sigma) \coloneqq g_n^{s,\sigma} \otimes g_n^{s,\sigma},
\]
are in $\mathscr C ^\infty \prth{\R^2, \L \prth{\R_y,\C^2} } $.
\end{theorem}

\noi Notably,
\[
    \mathrm T_\E(s,\sigma) = \sum_{n \in \Z} \lambda_n (s,\sigma) \Pi_n(s,\sigma).
\]

\begin{proof}
According to Remark~\ref{rem:spectrum}, we study $\mathscr T$.
Let $\lambda_\zeta \in \Sp(\mathscr T), \begin{pmatrix} f \\ g\end{pmatrix} \in E_{\lambda_\zeta}(\mathscr T)$ if and only if
\[
    \left\{ \begin{matrix} & & \mathfrak{a}^\dagger g & = & \left( \lambda_\zeta - \zeta \right) f \\ \mathfrak{a} f & - & 2 \zeta g & = & \left( \lambda_\zeta - \zeta \right) g \end{matrix}\right.
\]
If $\lambda_\zeta = \zeta$, then the injectivity of $\mathfrak{a}^\dagger$ on $\L\prth{\R}$ implies $g=0$ whence $\mathfrak{a}f = 0$ so $f = A \mathfrak{h}_0$ for some $A \in \C$.
So $\zeta \in \Sp(\mathscr T)$ and is a simple eigenvalue with $\displaystyle{ \begin{pmatrix} \mathfrak{h}_0 \\ 0 \end{pmatrix}}$ a normalized eigenfunction.
If $\lambda_\zeta \neq \zeta$, then
\[
    \left\{ \begin{matrix} f & = & \frac{1}{\lambda_\zeta - \zeta}\mathfrak{a}^\dagger g \\
    \mathfrak{a}\mathfrak{a}^\dagger g & = & \left( \lambda_\zeta^2 - \zeta^2 \right) g \end{matrix}\right.
\]
so, according to Proposition~\ref{pro:creation_annihilation}, we deduce that there exists $n \in \N_{> 0}, \lambda_\zeta^2 - \zeta^2 = 2n$ so
\[
    \lambda_\zeta \in \left\{ \sgn(n) \sqrt{\zeta^2 + 2|n|} \ | \ n \in \Z_{\neq 0} \right\},
\]
and for some $B \in \C$,
\[
    \forall (\zeta, y) \in \R^2, \ \begin{pmatrix} f_n \\ g_n \end{pmatrix} (\zeta, y) \coloneqq B \begin{pmatrix} \frac{\sqrt{2|n|}}{\sgn(n) \sqrt{\zeta^2 + 2|n|} - \zeta} \mathfrak{h}_{|n|} (y) \\ \mathfrak{h}_{|n|-1}(y) \end{pmatrix}.
\]
Then
\[
    \Sp \Big( \mathrm T_\E(s,\sigma) \Big) = \left\{ \lambda_0(\sigma) \right\} \bigcup \left\{ \lambda_n(s,\sigma) \ | \ n \in \Z_{\neq 0} \right\}.
\]
The associated projectors $\prth{\Pi_n}_{n \in \Z}$ are smooth because $r$ is smooth ($\m$ is in $\mathscr C ^\infty\prth{\R^2,\R}$) and strictly non-negative (Assumption~\ref{ass:transversality}); for all $n \in \Z_{\neq 0}$, $\lambda_n$ and $\lambda_n + \sigma$ is smooth and never vanishes.
\end{proof}

\section*{Acknowledgments}

The author is very grateful to his advisors Clotilde Fermanian-Kammerer and Nicolas Raymond for many insightful discussions and valuable suggestions.
The author acknowledges the support of the Région Pays de la Loire via the Connect Talent Project HiFrAn 2022 07750, and from the France 2030 program, Centre Henri Lebesgue ANR-11-LABX-0020-01.

\bibliographystyle{plain}
\bibliography{biblio}

\bigskip

(\'E. Vacelet) \textbf{Univ Angers}, CNRS, \textbf{LAREMA}, SFR MATHSTIC, F-49000 Angers, France,

\textit{Email address :} \textrm{eric.vacelet@univ-angers.fr}

\end{document}